\pdfoutput=1
\documentclass[11pt,a4paper]{article}
\usepackage[T1]{fontenc}
\usepackage[utf8]{inputenc}
\usepackage[english]{babel}
\usepackage{lmodern}
\usepackage[margin=27mm]{geometry}
\usepackage{microtype,amsmath,amssymb,amsthm,mathtools}
\usepackage{enumitem,booktabs,array,tabularx}
\usepackage[unicode,hidelinks]{hyperref}
\hypersetup{pdftitle={A negative answer to the final question in Les petits cailloux},pdfauthor={Mihai Prunescu}}
\allowdisplaybreaks
\setlist{itemsep=3pt,topsep=5pt}
\newtheorem{theorem}{Theorem}[section]
\newtheorem{proposition}[theorem]{Proposition}
\newtheorem{lemma}[theorem]{Lemma}
\newtheorem{corollary}[theorem]{Corollary}
\theoremstyle{definition}
\newtheorem{definition}[theorem]{Definition}
\newtheorem{question}[theorem]{Question}
\theoremstyle{remark}
\newtheorem{remark}[theorem]{Remark}
\newcommand{\A}{\mathcal A}
\newcommand{\M}{\mathcal M}
\newcommand{\Sstr}{\mathcal S}
\newcommand{\Pclass}{\mathsf P}
\newcommand{\NPclass}{\mathsf{NP}}
\newcommand{\NBPclass}{\mathsf{NBP}}
\newcommand{\coNBPclass}{\mathsf{coNBP}}
\newcommand{\PHclass}{\mathsf{PH}}
\newcommand{\Trace}{\operatorname{Trace}}
\newcommand{\EXPclass}{\mathsf{EXP}}
\newcommand{\Ppar}{\mathsf P_{\mathrm{par}}}
\newcommand{\NPpar}{\mathsf{NP}_{\mathrm{par}}}
\newcommand{\DTIME}{\operatorname{DTIME}}

\newcommand{\SAT}{\operatorname{SAT}}
\newcommand{\Bad}{\operatorname{Bad}}
\newcommand{\Step}{\operatorname{Step}}
\newcommand{\Acc}{\operatorname{Acc}}
\newcommand{\Init}{\operatorname{Init}}
\newcommand{\Addr}{\operatorname{Addr}}
\newcommand{\pack}{\operatorname{pack}}
\newcommand{\nil}{\mathsf{nil}}
\newcommand{\push}{\operatorname{push}}
\newcommand{\blank}{\mathsf{blank}}
\newcommand{\data}{\mathsf{data}}
\newcommand{\bits}{\{0,1\}}

\title{A negative answer to the final question in \emph{Les petits cailloux}}
\author{Mihai Prunescu\thanks{Research Center for Logic, Optimization and Security (LOS), Faculty of Mathematics and Computer Science, University of Bucharest, Academiei 14, Bucharest (RO-010014), Romania.}\thanks{Simion Stoilow Institute of Mathematics of the Romanian Academy, Research unit 5, P. O. Box 1-764, Bucharest (RO-014700), Romania. E-mail: \texttt{mihai.prunescu@imar.ro}.}}
\date{}

\begin{document}
\maketitle

\begin{abstract}
We consider polynomial-time computation over arbitrary structures, with unit cost for the primitive operations and relations and with nondeterministic witnesses consisting of arbitrary structure elements. We prove that every finite-signature structure with two distinct available elements and containing a binary operation $*$ and total unary operations $L,R$ such that $L(x*y)=x$ and $R(x*y)=y$ satisfies $\Pclass\ne\NPclass$. The proof encodes an exponentially long computation in one element. An incorrect transition has a polynomial-size binary address, and hence the existence of an error is an $\NBPclass$ problem. Existence of a correct accepting trace belongs to the mixed second level $\Sigma_2^{S,B}$, with an arbitrary-element existential block followed by a binary universal block. We obtain $\EXPclass\subseteq\Sigma_2^{S,B}\subseteq\PHclass$. A hypothetical equality $\Pclass=\NPclass$ would collapse this hierarchy to $\Pclass$ and imply $\Pclass=\EXPclass$, contradicting a direct deterministic diagonalization over the same structure. A simulation by function graphs extends the obstruction to every jointly injective binary primitive under the same assumptions on equality and available elements. Consequently no unary predicate on Malcev's free binary-tree algebra yields $\Pclass=\NPclass$, with or without the projections. Since Cantor pairing is first-order definable in standard arithmetic, every finite expansion of $(\mathbb N;+,\cdot,0,1)$ also satisfies the separation. In the elementary arithmetic signature consisting of addition, remainder, and base-two exponentiation, Cantor pairing and its projections are fixed term functions, so the trace inclusion holds there unconditionally. We compare the arithmetic truth-predicate obstruction with Tarski's theorem and, under elementary syntax coding, obtain a quantifier-free diagonal obstruction in the elementary arithmetic language. In contrast, starting from the author's unary example, we construct five infinite structures with binary operations and $\Pclass=\NPclass$. These include an algebra with one binary operation and one constant, an expansion with an associative rectangular-band operation, and an algebra with one binary operation whose every row and column has infinite image.


\end{abstract}

\section{Introduction}
\label{sec:introduction}

The model of computation introduced by Blum, Shub and Smale permits exact operations on real numbers at unit cost and gives structure-dependent versions of the usual complexity classes \cite{BSS}. The monograph of Blum, Cucker, Shub and Smale develops this approach systematically \cite{BCSS}. Goode extended the foundational study of computation from rings to arbitrary structures \cite{Goode1994}. Poizat's \emph{Les petits cailloux} places it in a model-theoretic setting: the computational universe is a first-order structure, the available elementary operations and tests are its functions and relations, and the length of an input is the number of its components \cite{Poizat1995}. One may therefore ask whether a particular infinite structure satisfies $\Pclass=\NPclass$ as a structure-relative question, distinct from the corresponding question for finite-alphabet Turing machines.

The choice of primitives matters. Meer studied restricted real machines and obtained separation results for suitable operation sets \cite{Meer1992,Meer1993}; Koiran analyzed computation with addition and order \cite{Koiran1994}. Fournier and Koiran subsequently compared lower-bound questions for additive real computation with their classical finite-alphabet counterparts, distinguishing the parameter-free and arbitrary-parameter settings \cite{FournierKoiran2000}. These works show why a statement about $\Pclass$ and $\NPclass$ over one signature cannot be transferred to another merely because the underlying set is unchanged. They also motivate the explicit treatment of constants and permitted instructions below.

In the general model-theoretic setting, the question whether an infinite structure of finite signature can satisfy $\Pclass=\NPclass$ has a positive answer. Hemmerling constructed suitable expansions of structures on binary words \cite{Hemmerling2005}, and Ga\ss ner gave related constructions \cite{Gassner2006}. The author constructed a structure with three unary functions and a unary predicate in which existential satisfaction with arbitrary parameters is decidable in polynomial time \cite{PrunescuJSL}. The connection between this property and fast elimination of quantifiers is explained in \cite{PrunescuElim}. The latter condition concerns a uniform decision procedure, or equivalently appropriate circuits; it does not require an explicitly written equivalent quantifier-free formula of polynomial length.

The construction in \cite{PrunescuJSL} develops Poizat's use of a truth predicate on a unary successor structure. In Poizat's earlier construction, an existential formula with one free variable is represented by a term whose truth is read by the predicate \cite{Poizat2000}. Unary terms cannot combine an arbitrary tuple of independent parameters into one element. The author's construction overcomes this limitation by computing finite local information about the parameter tuple and consulting a predicate that encodes suitable sentences. We recall this existence theorem in Section~\ref{sec:positive}, without repeating the proof.

There remains a distinct proposal: use the binary operation of Malcev's free tree algebra to combine parameters, and expand that algebra by a truth predicate. We consider the following precise version of Poizat's question, in the framework of \cite{Poizat1995,Poizat2000}; the specific Malcev proposal is recalled in \cite[p.~465]{PrunescuElim}.

\begin{question}[Les petits cailloux]
\label{question:poizat}
Let $\M$ be the algebra of finite binary trees, with the operation $x*y$ forming a new root with left subtree $x$ and right subtree $y$, and with two distinct named elements $0,1$. Does there exist a unary predicate $V$ such that the resulting structure satisfies $\Pclass=\NPclass$ for unit-cost computation?
\end{question}

We answer this question negatively. The result holds whether the left and right projections are primitive operations or are omitted. More generally, a finite-signature structure with two distinct available elements and a jointly injective binary primitive cannot satisfy $\Pclass=\NPclass$. The obstruction persists after adding any finitely many further operations or relations. It is therefore independent of the particular method proposed for constructing $V$. Our second running example is standard arithmetic
$\mathcal N=(\mathbb N;+,\cdot,0,1)$. Cantor pairing and its two
projections are first-order definable in this original signature.
A simulation by their graphs transfers the obstruction to every finite
expansion of $\mathcal N$, in particular to
$(\mathbb N;+,\cdot,0,1,V)$ for every unary predicate $V$.
Section~\ref{subsec:arithmetic-truth} compares the resulting restriction
on existential truth predicates with Tarski's undefinability theorem.
Section~\ref{subsubsec:elementary-arithmetic} then treats the signature
consisting of addition, remainder, and base-two exponentiation. In that
language the pairing and projections are fixed term functions, giving
an unconditional trace inclusion and a quantifier-free diagonal
obstruction under elementary syntax coding.

The reason is computational. With a pairing function and its projections, one structure element can encode any prescribed finite list of configurations, including an exponentially long list. An individual entry is accessible by a binary address of polynomial length. Thus a nondeterministic polynomial-time machine can find an error in the proposed computation. If $\Pclass=\NPclass$, the absence of errors can be checked deterministically in polynomial time. Guessing the whole trace as one element then makes every deterministic exponential-time computation an $\NPclass$ computation, and the assumed collapse makes it polynomial-time. A universal-machine diagonalization excludes this conclusion. Section~\ref{subsec:polynomial-hierarchy} defines the full and mixed polynomial hierarchies, and Lemma~\ref{lem:trace} locates the trace condition precisely in $\Sigma_2^{S,B}$.

The local checking of a computation is in the tradition of Cook's satisfiability construction \cite{Cook1971}. Its existential trace and universally checked positions also have the form familiar from alternating computation \cite{CKS1981}. The extra feature used here is specific to the structure model: the entire trace is a single permissible nondeterministic witness, however large its external representation. Related difficulties concerning compressed oracle encodings occur in Ga\ss ner's work \cite{Gassner2009}. Ga\ss ner also observes that, in her earlier construction on binary trees, equality of trees is not available as a one-step test \cite[p.~1204]{Gassner2009}. Equality is a primitive test throughout the present paper. We give the full argument, including the structure-relative hierarchy proof, so that no transfer of a classical theorem is left implicit.

The negative result does not exclude all binary operations. Starting from the unary structure $\A$ of \cite{PrunescuJSL}, we construct structures $\mathcal B,\mathcal C,\mathcal D,\mathcal E,\mathcal F$ with $\Pclass=\NPclass$. One of them has just one binary operation and one constant, with no primitive predicate. Another has an associative idempotent binary operation, and a further example has one binary operation with infinite image in every row and every column. Each construction is stated as a corollary with a separate simulation proof. 

\section{Computation and satisfaction over a structure}
\label{sec:computation}

\subsection{The free Malcev algebra and standard arithmetic}
\label{subsec:running-examples}

\begin{definition}[The Malcev algebra freely generated by a set]
\label{def:malcev}
Let $G$ be a nonempty set. The universe $T_G$ consists of all finite
rooted ordered full binary trees whose leaves are labelled by elements
of $G$. Equivalently, it is the set of formal terms generated by the rules
\[
 g\in T_G\quad(g\in G),\qquad
 u,v\in T_G\ \Longrightarrow\ (u*v)\in T_G.
\]
Different generators denote different elements, no generator is a
product, and
\[
 u*v=u'*v'\quad\Longleftrightarrow\quad u=u'\ \text{and}\ v=v'.
\]
The operation $*$ forms a new root with the two given trees as its
ordered subtrees. The algebra $(T_G;*)$ is the \emph{free Malcev algebra}
on $G$, also called the free magma on $G$: every map from $G$ to an
algebra with one binary operation extends uniquely to a homomorphism
from $T_G$, by recursion on finite terms. This is the case of a single
binary operation among the absolutely free algebras studied by Malcev
\cite[p.~730]{Malcev1962}. Here ``Malcev algebra'' refers to this free
tree algebra, not to the anticommutative algebras bearing the same name
in another part of algebra.
\end{definition}

Choose two distinct named elements $0,1\in T_G$. If $G$ has two chosen
generators they may be used; if $G=\{g\}$, the elements $g$ and $g*g$
may be used. The names need not denote leaves. Although $G$ may be
infinite or uncountable, we do not put a constant symbol for every
generator into the signature. All computational signatures below remain
finite. Define total projections by
\[
 L(u*v)=u,\qquad R(u*v)=v,\qquad L(g)=R(g)=g\quad(g\in G).
\]
We use the notation
\[
 \M=(T_G;*,0,1),\qquad \M^+=(T_G;*,L,R,0,1),
\]
and, for any unary predicate $V\subseteq T_G$,
\[
 \M_V=(T_G;*,V,0,1),\qquad
 \M_V^+=(T_G;*,L,R,V,0,1).
\]
Thus $*$ is a pairing function and $L,R$ are its left and right
projections. These names will be used consistently throughout the paper.

The second running example is \emph{standard arithmetic}
\[
 \mathcal N=(\mathbb N;+,\cdot,0,1),\qquad
 \mathbb N=\{0,1,2,\ldots\}.
\]
For $V\subseteq\mathbb N$, put
$\mathcal N_V=(\mathbb N;+,\cdot,0,1,V)$.
The pairing and projections used in the argument are definable functions
of this structure, rather than assumed primitives of its original
signature. Cantor's pairing function is
\[
 c(x,y)=\frac{(x+y)(x+y+1)}{2}+y.
\]
For each $n\in\mathbb N$, let $s$ be the unique nonnegative integer with
\[
 \frac{s(s+1)}2\leq n<\frac{(s+1)(s+2)}2,
\]
put $r=n-s(s+1)/2$, and define $L_c(n)=s-r$ and $R_c(n)=r$.
The interval has $s+1$ integers, so $0\leq r\leq s$. Consequently
\[
 L_c(c(x,y))=x,\qquad R_c(c(x,y))=y,\qquad
 c(L_c(n),R_c(n))=n.
\]
Indeed, $c(x,y)$ belongs to the displayed interval with $s=x+y$
and $r=y$; conversely, the formula reconstructs every $n$. This proves
that $c:\mathbb N^2\to\mathbb N$ is a bijection and that its two
coordinate maps are total.

For definability in the original arithmetic signature, use the
quantifier-free graph formula
\begin{equation}
 \Gamma_c(x,y,z)\quad\Longleftrightarrow\quad
 z+z=(x+y)\cdot(x+y+1)+(y+y).
 \label{eq:cantor-graph}
\end{equation}
In $\mathcal N$, this is equivalent to $z=c(x,y)$. The projection graphs
are defined by the existential arithmetic formulas
\[
 \Gamma_{L_c}(n,x)\ \Longleftrightarrow\ \exists y\,\Gamma_c(x,y,n),
 \qquad
 \Gamma_{R_c}(n,y)\ \Longleftrightarrow\ \exists x\,\Gamma_c(x,y,n).
\]
Bijectivity proves that these formulas define total single-valued
functions. No division or order symbol occurs in the graph definitions.
These same definitions remain valid in every expansion by additional
predicates or functions.

When the argument calls for these functions as primitive instructions,
we use only the \emph{auxiliary definitional expansion}
\[
 \mathcal N^{\mathrm{pair}}
   =(\mathbb N;+,\cdot,c,L_c,R_c,0,1),
 \qquad
 \mathcal N_V^{\mathrm{pair}}=(\mathcal N^{\mathrm{pair}},V).
\]
First-order definability alone does not supply a deterministic
unit-cost algorithm for producing a function's value. The graph of $c$
has a bounded-time test, and the projection graphs belong to
$\NPclass$ in every arithmetic expansion. Under a hypothetical
$\Pclass=\NPclass$, Lemma~\ref{lem:graphs} permits adjoining these
functions while preserving the collapse. Corollary~\ref{cor:arithmetic}
then transfers the pairing obstruction to the original arithmetic
signature. The running example remains $\mathcal N$.
Section~\ref{subsubsec:elementary-arithmetic} considers an elementary
arithmetic signature in which $c,L_c,R_c$ are fixed term functions.

\subsection{Machines, input length, and complexity classes}

Let $\mathcal L$ be a fixed finite signature of total functions, relations, and constants, and let $\Sstr$ be an $\mathcal L$-structure with universe $S$. Equality is always available as a test. We assume that two distinct elements, denoted by $0,1$, are named or obtainable by fixed ground terms. Replacing such ground terms by named constants changes running time by at most a fixed initialization cost.

A machine over $\Sstr$ has a finite program, finitely many states, and finitely many tapes. A data cell contains one element of $S$. The elementary instructions are copying a cell, writing a named constant, moving a tape head by one position, applying one primitive function, testing equality or a primitive relation, and testing whether a cell is blank. Each instruction costs one time unit. A fixed finite control alphabet and a blank marker are permitted; they do not supply additional elements of $S$. This is the multitape formulation used in \cite[Section~2]{PrunescuElim}.

An input is a tuple $\bar x\in S^n$, and its size is $|\bar x|=n$. In particular, neither the height of a tree nor the bit length of an external representation of a structure element contributes to input length. A primitive function call costs one step, but evaluating a variable-length term requires executing its individual operation occurrences. Finite descriptions of formulas and programs will be given as explicit binary words, whose lengths are charged in full. Here a binary description means a word over the two-letter alphabet $\{0,1\}$ encoding the syntax, instruction names, and indices; it is not an external representation of the structure elements on which the program operates.

\begin{definition}
The class $\DTIME_{\Sstr}(t(n))$ consists of the languages $X\subseteq S^*=\bigsqcup_{n\geq0}S^n$ decided by deterministic machines over $\Sstr$ in time $O(t(n))$ on every length-$n$ input. Define
\[
 \Pclass(\Sstr)=\bigcup_{k\geq1}\DTIME_{\Sstr}((n+2)^k),
 \qquad
 \EXPclass(\Sstr)=\bigcup_{k\geq1}\DTIME_{\Sstr}(2^{(n+2)^k}).
\]
Nondeterministic machines may guess arbitrary elements of $S$. The analogous polynomial-time class is $\NPclass(\Sstr)$. Equivalently, $X\in\NPclass(\Sstr)$ if there are a deterministic polynomial-time verifier $D$ and a polynomial $q$ such that
\[
 \bar x\in X\quad\Longleftrightarrow\quad
 \exists\bar y\in S^{q(|\bar x|)}\;D(\bar x,\bar y)=1.
\]
The two tuples are supplied with their lengths, or on separate input tapes; either convention is polynomially equivalent to a self-delimiting concatenation.
\end{definition}

All time bounds are uniform in the values of the input elements. Padding witnesses by a named constant allows one to replace ``length at most $q(n)$'' by ``length exactly $q(n)$''. Binary branching can be simulated by guessing a value and verifying that it is $0$ or $1$. A witness for $\NPclass$ need not be a binary word: its entries can be arbitrary elements of $S$. We distinguish this from binary nondeterminism in the next subsection.

At first, a program may use only the constants of the fixed signature. This is the convention of \cite{PrunescuElim}. The versions allowing a finite tuple of additional constants in each program are denoted by $\Ppar(\Sstr)$ and $\NPpar(\Sstr)$ and are discussed in Section~\ref{subsec:constants}. Both the positive and the negative conclusions of this paper hold under either convention.

\subsection{Binary nondeterminism and polynomial hierarchies}
\label{subsec:polynomial-hierarchy}

The bounded-quantifier description of the polynomial hierarchy is
classical \cite{Stockmeyer1976}; alternation gives its corresponding
machine interpretation \cite{CKS1981}. Full and binary, also called
digital, polynomial hierarchies over arbitrary structures are studied
in \cite{BournezCuckerNauroisMarion2006}. We specify the quantifier
domains because unrestricted structure-element guesses and binary
guesses have different roles in the trace argument. The mixed-domain
notation below is introduced explicitly for this paper.

\begin{definition}[Binary nondeterminism]
The class $\NBPclass(\Sstr)$ consists of languages $X\subseteq S^*$ for
which there are a deterministic polynomial-time verifier $D$ over
$\Sstr$ and a polynomial $q$ such that, with $n=|\bar x|$,
\[
 \bar x\in X\quad\Longleftrightarrow\quad
 \exists\bar\varepsilon\in\bits^{q(n)}\;
 D(\bar x,\bar\varepsilon)=1.
\]
Thus the \emph{witness}, but not necessarily the input, is a binary word.
Write $\coNBPclass(\Sstr)$ for the complements of these languages, with
complement taken in $S^*$. We have
\[
 \Pclass(\Sstr)\subseteq\NBPclass(\Sstr)
 \subseteq\NPclass(\Sstr).
\]
The second inclusion follows by verifying that every guessed entry is
either $0$ or $1$.
\end{definition}

\begin{definition}[Full and mixed polynomial hierarchies]
\label{def:polynomial-hierarchy}
Set $\Sigma_0^{\Pclass}(\Sstr)=\Pi_0^{\Pclass}(\Sstr)=\Pclass(\Sstr)$.
For a fixed $k\geq1$, a language $X\subseteq S^*$ belongs to
$\Sigma_k^{\Pclass}(\Sstr)$ if there are polynomial bounds
$p_1,\ldots,p_k$ and a deterministic polynomial-time predicate $D$ such
that, for $n=|\bar x|$,
\[
 \bar x\in X\quad\Longleftrightarrow\quad
 \exists\bar y_1\in S^{p_1(n)}\;
 \forall\bar y_2\in S^{p_2(n)}\;\cdots\;
 Q_k\bar y_k\in S^{p_k(n)}\;
 D(\bar x,\bar y_1,\ldots,\bar y_k)=1,
\]
where the blocks alternate and $Q_k$ is existential for odd $k$ and
universal for even $k$. The class $\Pi_k^{\Pclass}(\Sstr)$ is defined by
starting with a universal block, or equivalently by complementation.
In particular,
\[
 \Sigma_1^{\Pclass}(\Sstr)=\NPclass(\Sstr),\qquad
 \Pi_1^{\Pclass}(\Sstr)=\operatorname{coNP}(\Sstr).
\]
The full polynomial hierarchy is
\[
 \PHclass(\Sstr)=\bigcup_{k\geq0}
 \bigl(\Sigma_k^{\Pclass}(\Sstr)\cup\Pi_k^{\Pclass}(\Sstr)\bigr).
\]
All tuple boundaries are explicit, all bounds concern numbers of cells,
and the running time of $D$ is polynomial in the total input length.
Polynomial witness bounds are taken to have nonnegative integer values
at every input length.

For a \emph{mixed} hierarchy, let $B=\bits$ and choose a domain
$\tau_i\in\{S,B\}$ for each block. Replacing $S^{p_i(n)}$ above by
$\tau_i^{p_i(n)}$ defines our notation
$\Sigma_k^{\tau_1,\ldots,\tau_k}(\Sstr)$, and starting universally
defines $\Pi_k^{\tau_1,\ldots,\tau_k}(\Sstr)$. Thus the all-$S$ classes
are the full classes, $\Sigma_1^B(\Sstr)=\NBPclass(\Sstr)$, and the
all-$B$ classes use only binary quantifier blocks. In particular, a
language $X$ belongs to $\Sigma_2^{S,B}(\Sstr)$ if and only if there
are fixed polynomials $p,q$ and a deterministic polynomial-time
predicate $D$ such that, for every input $\bar x$ of length $n$,
\begin{equation}
 \bar x\in X\quad\Longleftrightarrow\quad
 \exists\bar y\in S^{p(n)}\;
 \forall\bar\varepsilon\in B^{q(n)}\;
 D(\bar x,\bar y,\bar\varepsilon)=1.
 \label{eq:mixed-sigma-two}
\end{equation}
\end{definition}

Every mixed class at level $k$ is contained in the corresponding full
class. For an existential binary block, quantify over $S$ and require
that every entry be $0$ or $1$; for a universal binary block, use that
same condition as the antecedent of an implication. These guards have
polynomial-time tests and add no quantifier alternation. In particular,
\[
 \Sigma_2^{S,B}(\Sstr)\subseteq\Sigma_2^{\Pclass}(\Sstr)
 \subseteq\PHclass(\Sstr).
\]
The matrices in these definitions are polynomial-time predicates; we
are not requiring one fixed quantifier-free formula to express every
variable-length computation of such a matrix.

\begin{proposition}[Collapse of the full and mixed hierarchies]
\label{prop:ph-collapse}
If $\Pclass(\Sstr)=\NPclass(\Sstr)$, then
\[
 \PHclass(\Sstr)=\Pclass(\Sstr).
\]
Every fixed level of every mixed hierarchy also equals
$\Pclass(\Sstr)$. In particular,
$\Sigma_2^{S,B}(\Sstr)=\Pclass(\Sstr)$.
\end{proposition}
\begin{proof}
Eliminate the quantifier blocks from the inside out. An existential
block of polynomially many structure elements applied to a polynomial-time
predicate defines an $\NPclass(\Sstr)$ relation, and hence, under the
hypothesis, a $\Pclass(\Sstr)$ relation. For a universal block, complement
the predicate, eliminate the resulting existential block, and complement
the answer. Deterministic polynomial time is closed under complementation.
Previously quantified tuples remain part of the input, with their
boundaries and the original input length retained, so every witness
bound stays polynomial in the current input length. This proves by
induction that every fixed full level is contained in $\Pclass(\Sstr)$.
The reverse inclusions use predicates that ignore their quantified
variables. The mixed classes are contained in the corresponding full
classes, and the same argument with unused variables proves equality
for them. The number of blocks is fixed for each language; no uniform
polynomial bound for unbounded alternation is asserted.
\end{proof}

\subsection{The problem \texorpdfstring{$\SAT(\Sstr)$}{SAT(S)}}

\begin{definition}
An instance of $\SAT(\Sstr)$ consists of an explicit binary encoding $w_\psi$ of an existential $\mathcal L$-formula
\[
 \psi(\bar x)=\exists\bar y\,\theta(\bar x,\bar y),
\]
where $\theta$ is quantifier-free, together with a tuple $\bar a\in S^{|\bar x|}$. The instance is accepted exactly when $\Sstr\models\psi(\bar a)$. Its length is $|w_\psi|+|\bar a|$, up to a fixed choice of delimiters. Malformed encodings or incorrect parameter arities are rejected. The matrix may contain arbitrary Boolean combinations and negative occurrences of every relation symbol, including $V$ when it is present.
\end{definition}

\begin{proposition}[Structure-relative satisfiability]
\label{prop:sat}
The problem $\SAT(\Sstr)$ is $\NPclass(\Sstr)$-complete under polynomial-time many-one reductions over $\Sstr$. Consequently,
\[
 \Pclass(\Sstr)=\NPclass(\Sstr)
 \quad\Longleftrightarrow\quad
 \SAT(\Sstr)\in\Pclass(\Sstr).
\]
\end{proposition}
\begin{proof}
For membership, guess the values of the quantified variables and evaluate the explicit quantifier-free matrix. Variable names and the syntax tree can be processed with polynomial overhead; each operation or relation occurrence requires one primitive call after its arguments have been evaluated.

For hardness, fix a nondeterministic polynomial-time machine and a polynomial clock $T(n)\geq n$ bounding every branch. During $T(n)$ steps, every tape uses at most $O(T(n))$ relevant cells. Describe a computation tableau with one row for each time and one column for each relevant tape cell. Introduce variables for the data values, blank flags, head positions, states, and nondeterministic choices. There are polynomially many variables and local conditions. Finite control values are represented by fixed tuples of binary variables, each constrained by $b=0\lor b=1$. A function instruction is imposed by an equation $z=f(u_1,\ldots,u_r)$. A relation test with a binary result $b$ is imposed by
\[
 (b=1\land R(\bar u))\lor(b=0\land\neg R(\bar u)).
\]
Copying, unchanged tape cells, tape movements, and conditional instructions have fixed-size local descriptions. Require the initial row to contain the given input tuple as its free parameters and the final row to be accepting, after making accepting states absorbing. Existentially quantify all internal data and control variables. The resulting formula is satisfiable precisely when the original machine has an accepting branch.

The syntax of this formula is produced in polynomial time, independently of the values of the input parameters. Every intermediate value has its own variable, so shared computation is not expanded into an exponentially long term. This is the Cook construction \cite{Cook1971} over the primitive operations of $\Sstr$; see \cite[Theorems~2--3]{PrunescuElim} for its structure-relative formulation.
\end{proof}

\begin{remark}
The witness variables in this proof need not be binary. In particular, the proposition is about arbitrary-element nondeterminism. We also make no general assertion that $\NPclass(\Sstr)\subseteq\EXPclass(\Sstr)$: over an arbitrary infinite structure there need not be an effective bounded enumeration of all possible witness values.
\end{remark}

\subsection{The deterministic separation \texorpdfstring{$\Pclass\ne\EXPclass$}{P differs from EXP}}
\label{subsec:hierarchy}

The following argument is the padded universal-machine diagonalization underlying the deterministic time hierarchy of Hartmanis and Stearns \cite{HartmanisStearns1965}. We supply the proof over $\Sstr$, because its primitive operations and relations may have no effective presentation as ordinary Turing-computable functions or predicates.

\begin{lemma}[Universal simulation]
\label{lem:universal}
For a fixed finite signature there is a universal machine over $\Sstr$ which, given a binary program description $e$, an input $\bar x$, and a bound $t$ on the number of simulated steps, simulates those steps in time
\[
 (|e|+|\bar x|+t+2)^{O(1)}.
\]
The exponent is independent of $e$, $\bar x$, and $t$.
\end{lemma}
\begin{proof}
Use a finite binary instruction format specifying states, tape indices, and one of the finitely many operation or relation symbols. The universal machine stores the simulated tapes, their head markers, and the current state. In $t$ steps the number of represented tape positions is bounded by a polynomial in $|e|+|\bar x|+t$. A full scan of this representation locates the current cells of the tapes addressed by an instruction. A scan of the finite description finds the instruction and its possible successor states.

The number of arguments of a primitive is bounded by the maximum arity in the fixed signature. The universal machine copies those arguments into its working cells and applies that same primitive operation or test. Dispatch among the finitely many primitive names is part of its finite program. Updating the tape representation and its markers requires at most another polynomial number of moves. This deliberately inefficient implementation has polynomial total overhead, with a fixed exponent. It interprets primitive names, not external representations of their values; no ordinary algorithm for a primitive's interpretation is required.
\end{proof}

\begin{theorem}[Deterministic separation]
\label{thm:hierarchy}
For every finite-signature structure with two distinct available elements,
\[
 \Pclass(\Sstr)\subsetneq\EXPclass(\Sstr).
\]
\end{theorem}
In particular, the conclusion applies to the free Malcev algebra,
standard arithmetic $\mathcal N$, and all their finite expansions,
including the structures with primitive pairing and projections used
in Section~\ref{sec:pairing}.
\begin{proof}
Choose an effective binary coding of finite programs. A padded description has the form
\[
 w=1^{|e|}0\,e\,0^r,\qquad r\geq0,
\]
where $e$ is a valid program code. Thus every fixed program has descriptions of arbitrarily large length. On a nonbinary or malformed input, the diagonal machine rejects. On a valid $w$ of length $n$, it simulates the program described by $e$ on the entire input $w$ for at most $2^n$ simulated steps. If the simulated program accepts within the bound, it rejects; if that program rejects within the bound, it accepts; on timeout it rejects.

An $(n+1)$-bit clock counts the simulated steps. Lemma~\ref{lem:universal}, together with the polynomial overhead for the clock, gives total running time $2^{O(n)}$ times a polynomial factor. The diagonal language $H$ therefore belongs to $\EXPclass(\Sstr)$.

Suppose that a machine $D$ decides $H$ in polynomial time. Let $e_D$ be its fixed program code and choose $r$ so large that the running time of $D$ on length
$n=2|e_D|+1+r$ is less than $2^n$. On
$w=1^{|e_D|}0e_D0^r$, the diagonal machine then observes $D$ halt and gives the opposite answer. This contradicts the assumption that $D$ decides $H$. Hence $H\notin\Pclass(\Sstr)$. The non-strict inclusion follows directly from the definitions.
\end{proof}

\subsection{Finitely many program constants}
\label{subsec:constants}

\begin{lemma}
\label{lem:constants}
If $\Pclass(\Sstr)=\NPclass(\Sstr)$, then the same equality holds after naming any finite tuple $\bar c\in S^m$. Moreover,
\[
 \Ppar(\Sstr)=\NPpar(\Sstr)
 \quad\Longleftrightarrow\quad
 \text{for some finite $\bar c$, }\Pclass(\Sstr,\bar c)=\NPclass(\Sstr,\bar c).
\]
\end{lemma}
\begin{proof}
For the first assertion, use Proposition~\ref{prop:sat}. In an existential formula of the expanded language replace each new constant symbol by a fresh free variable and supply its actual value as an additional parameter cell. This is polynomial syntactic preprocessing, and the given solver for $\SAT(\Sstr)$ decides the resulting instance. The number of added parameter cells is fixed.

If the parameter-allowing classes collapse, the base-language problem $\SAT(\Sstr)$ has a deterministic polynomial-time solver using some finite tuple $\bar c$ of program constants. Name that tuple. The same substitution reduces $\SAT(\Sstr,\bar c)$ to the base-language problem, giving a fixed-signature collapse. Conversely, suppose a finite expansion by $\bar c$ has a collapse. A given parameter-allowing nondeterministic machine uses some further finite tuple $\bar d$. The first assertion gives a collapse after naming both tuples, and therefore a deterministic polynomial-time decider with finitely many program constants. This is the finite-constant reduction discussed in \cite[p.~461]{PrunescuElim}.
\end{proof}

Thus any obstruction proved for every finite constant expansion also excludes the parameter-allowing collapse. Any positive fixed-signature example gives a positive example under the parameter-allowing convention as well.

\section{Pairing, computation traces, and the problem \texorpdfstring{$\Bad(A)$}{Bad(A)}}
\label{sec:pairing}

Throughout this section $\Sstr$ is a finite-signature structure containing distinct available elements $0,1$, a binary operation $*$, and total unary operations $L,R$ satisfying
\begin{equation}
 L(u*v)=u,\qquad R(u*v)=v\qquad(u,v\in S).
 \label{eq:pairing}
\end{equation}
We call $*$ a pairing function and $L,R$ its projections. These identities imply that $*$ is jointly injective. No other identity, in particular no identity $L(z)*R(z)=z$ for arbitrary $z$, is assumed. Additional finitely many operations and relations are allowed. They will be used in the transition checks exactly as they are interpreted in $\Sstr$.

The free Malcev algebra with its projections and the auxiliary
definitional expansion $\mathcal N^{\mathrm{pair}}$ of standard
arithmetic satisfy \eqref{eq:pairing}. In the latter, $*$ is interpreted
as $c$, and $L,R$ as $L_c,R_c$. It is not a free tree algebra:
$c(0,0)=0$ and $L_c(0)=R_c(0)=0$, so repeated projections may cycle.
The argument uses only \eqref{eq:pairing}, totality, and equality; it
requires neither freeness nor well-foundedness of arbitrary elements.
In this section the projections are primitive instructions. We do not
infer the same cost bounds for their evaluation in the original
$\{+,\cdot,0,1\}$-reduct from definability alone; the transfer to that
reduct is proved in Corollary~\ref{cor:arithmetic}. The unconditional
simulation by terms in elementary arithmetic is given in
Section~\ref{subsubsec:elementary-arithmetic}.

\subsection{Arrays represented by one element}

For $b_1\cdots b_d\in\bits^d$, with the most significant bit first, put
\[
 \Addr_d(T;b_1\cdots b_d)
 =P_{b_d}(\cdots P_{b_2}(P_{b_1}(T))\cdots),
 \qquad P_0=L,\quad P_1=R.
\]
For $0\leq i<2^d$, we also write $\Addr_d(T;i)$, using the length-$d$ binary expansion of $i$, including leading zeroes. The integer notation refers to ordinary binary control; it introduces no arithmetic operation on $S$.

\begin{lemma}[Packing and access]
\label{lem:array}
Given any list $C_0,\ldots,C_{2^d-1}$ of elements of $S$, there is an element $T\in S$ such that $\Addr_d(T;i)=C_i$ for every $i$. An entry can be accessed in $O(d)$ steps when the address is supplied in $d$ binary cells. Conversely, every $T\in S$ determines such a list of extracted values, without any shape assumption on $T$.
\end{lemma}
\begin{proof}
For $d=0$ take $T=C_0$. Inductively pack the first and second halves of the list into $T_0,T_1$ and set $T=T_0*T_1$. The first address bit selects the correct half by \eqref{eq:pairing}, and the induction hypothesis handles the remaining bits. The access algorithm reads one bit and applies its indicated projection at each step. Finally, $L$ and $R$ are total, so every sequence of $d$ projection calls has a value on every $T$.
\end{proof}

The existence part of the lemma does not give a polynomial-time construction of an arbitrary exponential array. The trace algorithm below guesses its root as one arbitrary structure element. This distinction is essential.

\subsection{Canonical configurations and transitions}
\label{subsec:configurations}

Fix a deterministic machine $A$ over the full signature of $\Sstr$. Here $A$ denotes a program; the calligraphic letter $\A$ in Section~\ref{sec:positive} will denote the unary structure. Encode every auxiliary finite control symbol by a fixed-size block over $\bits$, using extra tapes if needed. Since $A$ is fixed, this normalization has constant overhead and lets us assume that every nonblank tape cell contains an element of $S$. We describe a configuration encoding for which one transition has a bounded-time test, independently of the sizes of the encoded tapes.

For a fixed number $m$ of fields, use the tuple code
\[
 \pack_m(z_1,\ldots,z_m)
   =z_1*(z_2*(\cdots*(z_m*0)\cdots)).
\]
Its fields are recovered by fixed projection terms, and the tuple map is injective. Encode the finitely many states by distinct fixed-length binary tuples, packed in this way. Distinctness follows from \eqref{eq:pairing}, without any freeness assumption.

Use tagged stacks, writing $\push(a,s)$ for the result of placing
the element $a$ at the top of the stack $s$:
\begin{equation}
 \nil=0*0,\qquad \push(a,s)=1*(a*s).
 \label{eq:stacks}
\end{equation}
For a nonempty valid stack $s$, its top and remaining stack are
\[
 \operatorname{head}(s)=L(R(s)),\qquad
 \operatorname{tail}(s)=R(R(s)).
\]
The equality $\nil=\push(a,s)$ is impossible because their left projections are $0$ and $1$. This works even if the distinguished elements themselves belong to the image of $*$. In particular, it works in $\mathcal N^{\mathrm{pair}}$, where every natural number belongs to the image of Cantor pairing. Finite stacks remain
unambiguously decodable because the empty and nonempty tags are distinct;
no global condition on the shapes of arbitrary elements is used.

Separately tag tape cells by
\[
 \blank=0*0,\qquad \data(a)=1*a.
\]
A blank is distinct from every tagged data cell, including a cell whose data value is $0$, $1$, or one of the encoding tags. A tape is represented by three fields $(\ell,c,r)$: the stack to the left of the head, the current tagged cell, and the stack to its right. In either stack the nearest cell comes first. An empty stack represents an unvisited blank half-tape.

Popping an empty stack returns $(\blank,\nil)$; popping a nonempty
valid stack returns its head and tail. For an arbitrary $s\ne\nil$,
use the same projection formulas to define $\operatorname{pop}(s)$.
This fixes a total convention on malformed stacks without asserting
their validity. A move to the right replaces $(\ell,c,r)$ by
\[
 \bigl(\push(c,\ell),c',r'\bigr),
 \qquad(c',r')=\operatorname{pop}(r).
\]
A move to the left is symmetric. A data write changes only the current cell, and primitive tests are applied to the decoded current data values. Operations on blank operands may be normalized to a fixed rejecting transition, after an explicit blank test. This only fixes the semantics on cases where a data instruction has no operand and does not restrict the computations under consideration.

A \emph{valid configuration} is the canonical packed tuple of a state and the tape fields of $A$, with finite tagged stacks as just described. A represented tape may retain extra blank cells that have been visited; its encoding is part of the configuration. The successor update is fixed by the stack rules, so there is a unique encoded successor even when two different encodings would describe the same extensional infinite blank tape.

\begin{lemma}[Exact local transition test]
\label{lem:transition}
For the fixed machine $A$ there are constant-time tests $\Step_A(C,C')$ and $\Acc_A(C)$ with the following properties:
\begin{enumerate}
\item if $C$ is valid, exactly one $C'$ satisfies $\Step_A(C,C')$, namely the canonical encoded successor of $C$;
\item that successor is valid;
\item on valid configurations, $\Acc_A$ is precisely the accepting-state test.
\end{enumerate}
The initial configuration $\Init_A(\bar x)$ is constructible in $O(|\bar x|+1)$ steps. The accepting and rejecting halting states can be made absorbing.
\end{lemma}
\begin{proof}
The number of fields in a configuration depends only on $A$. Extract
them by the fixed projection terms. If the decoded state is not one
of the machine's state codes, the transition test returns false.
Otherwise, perform the instruction's actual primitive operation or
relation test. The formulas for stack movement, copying, and writing
determine all updated fields using a bounded number of primitive calls.

For each instruction and branch, require $C'$ to equal the \emph{whole canonical tuple} obtained by packing those updated fields. This equality is a single structure equality test after a bounded number of constructions. Taking the finite disjunction of the instruction cases defines $\Step_A$. On a valid $C$, determinism chooses exactly one case; the stack formulas preserve validity and give its unique encoded successor. Testing the decoded state gives $\Acc_A$. Self-loops at halting states are implemented by $C'=C$ on those states.

To initialize the input tape, tag the input elements and build its right stack from right to left. Its left stack is empty, and the remaining tapes are blank with empty stacks. This requires a linear number of pairing calls. The other fields take constant time.

No global validity test on arbitrary stacks is needed. The assertions concern successors of valid configurations, and validity will propagate from the known initial configuration. The whole-target equality is nevertheless indispensable: checking only the projections of $C'$ could admit noncanonical elements, since \eqref{eq:pairing} does not imply $L(z)*R(z)=z$ for every $z$.
\end{proof}

\subsection{Errors in an exponentially long trace}

For $d\geq1$, a proposed initial configuration $I\in S$, and a proposed trace $T\in S$, put
\[
 C_i=\Addr_d(T;i)\qquad(0\leq i<2^d).
\]
Define the language $\Bad(A)$ on inputs $(1^d,I,T)$ by
\begin{equation}
\begin{split}
 (1^d,I,T)\in\Bad(A)\quad\Longleftrightarrow\quad &
 C_0\ne I\ \lor\ \neg\Acc_A(C_{2^d-1})\\
 &{}\lor\exists i<2^d-1\;\neg\Step_A(C_i,C_{i+1}).
\end{split}
\label{eq:bad}
\end{equation}
The last two input cells contain $I,T$, and the preceding cells form the unary word $1^d$. Inputs outside this format are, by convention, in $\Bad(A)$. Thus its input length is exactly $d+2$ on well-formed instances, independently of the represented size of $I$ and $T$. The machine $A$ is fixed, not part of this input.

\begin{lemma}[Finding an error]
\label{lem:bad}
For each fixed $A$,
\[
 \Bad(A)\in\NBPclass(\Sstr)\subseteq\NPclass(\Sstr).
\]
Its verifier uses only polynomially many binary choices and has a
polynomial running-time bound in $d+2$, uniform in $I,T$. Consequently
the no-error language is in $\coNBPclass(\Sstr)$.
\end{lemma}
\begin{proof}
First check the input format in linear time and accept immediately if it fails. Otherwise choose one of the three error types in \eqref{eq:bad} using a fixed number of binary branches.

For the initial error, apply $L$ exactly $d$ times to $T$, compare the result with $I$, and accept if they differ. For the final error, apply $R$ exactly $d$ times and accept if $\Acc_A$ fails. These branches take $O(d+1)$ primitive steps apart from ordinary tape management.

For a transition error, guess a word $b_1\cdots b_d$ and verify that each entry is $0$ or $1$. Reject the branch if all entries are $1$, since that word denotes the last position and has no successor in the array. Compute the binary successor word by carrying from the last bit; this uses $O(d)$ bit operations and requires no arithmetic on $S$. Copy $T$ into two working cells and follow the two address words to obtain $C_i$ and $C_{i+1}$. Apply the constant-time test of Lemma~\ref{lem:transition} and accept precisely when it fails.

Every branch has a polynomial bound in $d+2$. If an endpoint is wrong, its corresponding branch accepts. If a transition is wrong, the branch guessing its index accepts. Conversely, an accepting branch certifies one of the displayed errors. This proves exact membership in $\Bad(A)$.
\end{proof}

Define the accepted-trace language on well-formed inputs by
\[
 \Trace(A)=\{(1^d,I):d\geq1,\ \exists T\in S\;
                      (1^d,I,T)\notin\Bad(A)\},
\]
and reject inputs outside this format. Its input length is $d+1$.
Here $I$ is any element of $S$: validity of $I$ is not included as an
extra membership test. The following lemma gives the interpretation
when $I$ is a valid configuration.

\begin{lemma}[Exact traces and the mixed second level]
\label{lem:trace}
For each fixed $A$,
\[
 \Trace(A)\in\Sigma_2^{S,B}(\Sstr)
              \subseteq\Sigma_2^{\Pclass}(\Sstr).
\]
More precisely, there is a deterministic polynomial-time predicate $Q_A$
such that
\begin{equation}
 (1^d,I)\in\Trace(A)\quad\Longleftrightarrow\quad
 \exists T\in S\;\forall\bar\varepsilon\in\bits^d\;
 Q_A(1^d,I,T,\bar\varepsilon).
 \label{eq:trace-hierarchy}
\end{equation}
If $I$ is a valid configuration of $A$, then for every $d\geq1$,
\begin{equation}
 \exists T\in S\;\bigl((1^d,I,T)\notin\Bad(A)\bigr)
 \quad\Longleftrightarrow\quad
 A\text{ accepts from }I\text{ within }2^d-1\text{ transitions}.
 \label{eq:trace}
\end{equation}
\end{lemma}
\begin{proof}
For a well-formed input and a binary address
$\bar\varepsilon\in\bits^d$, let $Q_A$ assert
\[
 \begin{gathered}
 \Addr_d(T;0^d)=I,\qquad
 \Acc_A\bigl(\Addr_d(T;1^d)\bigr),\\
 \bar\varepsilon=1^d\quad\text{or}\quad
 \Step_A\bigl(\Addr_d(T;\bar\varepsilon),
              \Addr_d(T;\bar\varepsilon+1)\bigr).
 \end{gathered}
\]
The last test is performed only when $\bar\varepsilon\ne1^d$;
$\bar\varepsilon+1$ denotes the binary successor word computed by the
carry procedure in Lemma~\ref{lem:bad}. The endpoints and the two
addressed configurations are accessed using $O(d)$ projection calls,
and Lemma~\ref{lem:transition} supplies the local tests. Thus $Q_A$ is
decidable in polynomial time in $d+2$, uniformly in $I,T$. Its matrix
rejects malformed inputs. To use quantifier lengths defined for every
input length $n$, give the universal word length $n$ and, on a
well-formed input $n=d+1$, ignore its last bit.

Every possible transition position is represented by a binary address,
while the all-one address requires only the endpoint checks. Therefore
the universal assertion is equivalent to the absence of all three
errors in \eqref{eq:bad}. This proves \eqref{eq:trace-hierarchy} and
the mixed-class membership. The existential block contains just one
arbitrary structure element. We now establish the claimed acceptance
equivalence for valid $I$.

Suppose first that $T$ has no error. Its first extracted entry is $C_0=I$, hence is valid. We prove by induction on $i$ that $C_i$ is the actual canonical configuration reached after $i$ transitions. The assertion holds at $i=0$. If it holds at $i<2^d-1$, the absence of a transition error gives $\Step_A(C_i,C_{i+1})$. Lemma~\ref{lem:transition} makes $C_{i+1}$ the unique actual successor and makes it valid. This proves the induction. The last entry satisfies $\Acc_A$, so the machine accepts within the stated bound.

Conversely, suppose that $A$ accepts after $h\leq2^d-1$ transitions. List its actual configurations through time $h$, and repeat the accepting configuration until the list has exactly $2^d$ entries. These repetitions are legal because the accepting state is absorbing. Lemma~\ref{lem:array} supplies a root $T$ whose extracted entries are exactly this list. Its initial entry is $I$, its final entry is accepting, and all adjacent pairs pass the transition test. None of the alternatives in \eqref{eq:bad} holds.

The forward argument applies to every candidate $T$, including elements that do not have a recognizable full-tree shape. Only the extracted entries and their exact successor relation are used.
\end{proof}

\subsection{The obstruction}

\begin{proposition}[Exponential time in the mixed second level]
\label{prop:exp-mixed}
For every structure satisfying the assumptions of this section,
\[
 \EXPclass(\Sstr)\subseteq\Sigma_2^{S,B}(\Sstr)
                  \subseteq\Sigma_2^{\Pclass}(\Sstr).
\]
In the first inclusion a single arbitrary structure element suffices
for the existential block.
\end{proposition}
\begin{proof}
Let $X\in\EXPclass(\Sstr)$ be decided by a fixed deterministic machine
$A$ in time at most $2^{q(n)}$, after increasing the polynomial $q$ to
absorb constants and small input lengths. Set
\[
 d=q(n)+1,\qquad I=\Init_A(\bar x).
\]
The map $\bar x\mapsto(1^d,I)$ is computable in polynomial time and
has polynomial output length. The initial configuration is valid, and
\[
 2^{q(n)}\leq2^{q(n)+1}-1=2^d-1.
\]
As both accepting and rejecting halting states are absorbing,
Lemma~\ref{lem:trace} gives
\[
 \bar x\in X\quad\Longleftrightarrow\quad
 \exists T\in S\;\forall\bar\varepsilon\in\bits^{q(n)+1}\;
 Q_A\bigl(1^{q(n)+1},\Init_A(\bar x),T,\bar\varepsilon\bigr).
\]
Its deterministic matrix first performs the polynomial-time
preprocessing and then the test $Q_A$. This is the required mixed
second-level definition, without any collapse assumption.
\end{proof}

\begin{theorem}[Pairing obstruction]
\label{thm:pairing}
Every finite-signature structure $\Sstr$ satisfying \eqref{eq:pairing}, with equality and two distinct available elements, satisfies
\[
 \boxed{\Pclass(\Sstr)\ne\NPclass(\Sstr).}
\]
The conclusion remains valid after any finite expansion by additional total functions, relations, or constants, and under the convention allowing finitely many arbitrary program constants.
\end{theorem}
\begin{proof}
Assume $\Pclass(\Sstr)=\NPclass(\Sstr)$. Proposition~\ref{prop:ph-collapse}
and Proposition~\ref{prop:exp-mixed} give
\[
 \EXPclass(\Sstr)\subseteq\Sigma_2^{S,B}(\Sstr)
 \subseteq\PHclass(\Sstr)=\Pclass(\Sstr),
\]
contrary to Theorem~\ref{thm:hierarchy}. To spell out the two
nondeterministic steps underlying this hierarchy argument, let
$X\in\EXPclass(\Sstr)$ and fix a deterministic decider $A$ for $X$ whose running time on length-$n$ inputs is at most $2^{q(n)}$, for an integer-valued polynomial $q$. Increasing $q$ absorbs constant factors and the finitely many small input lengths.

By Lemma~\ref{lem:bad}, the assumed equality gives a deterministic polynomial-time decider for $\Bad(A)$. Complementing its output gives such a decider for the absence of errors. On input $\bar x\in S^n$, perform the following nondeterministic computation: construct
\[
 I=\Init_A(\bar x),\qquad d=q(n)+1,
\]
write $1^d$, guess one arbitrary element $T\in S$, and use the deterministic decider to check that $(1^d,I,T)\notin\Bad(A)$. All deterministic work is polynomial in $n$, since the tested instance has length $d+2$.

By Lemma~\ref{lem:trace}, this procedure accepts exactly when $A$ accepts $\bar x$ within $2^d-1$ transitions. Since
\[
 2^{q(n)}\leq 2^{q(n)+1}-1,
\]
and rejecting states are absorbing, this is exactly the language $X$. Therefore $X\in\NPclass(\Sstr)=\Pclass(\Sstr)$. As $X$ was arbitrary, $\EXPclass(\Sstr)\subseteq\Pclass(\Sstr)$, contradicting Theorem~\ref{thm:hierarchy}.

Any finite expansion still has a finite signature and the identities \eqref{eq:pairing}; its actual additional primitives can be included in both the universal simulation and the transition test. The same proof therefore applies to it. Finally, a parameter-allowing collapse would yield a collapse in some finite constant expansion by Lemma~\ref{lem:constants}, which has just been excluded.
\end{proof}

\begin{corollary}[No fast existential truth predicate]
\label{cor:truth}
No finite expansion of a structure satisfying \eqref{eq:pairing} admits a uniform polynomial-time algorithm for $\SAT$ in its full signature. In particular, no unary predicate $V$ in such an expansion can satisfy
\[
 V(C_\psi(\bar a))\quad\Longleftrightarrow\quad\psi(\bar a)
\]
for every existential formula of the expanded language and every parameter tuple, when $C_\psi(\bar a)$ is uniformly constructible in polynomial unit cost from the explicit formula and the tuple.
\end{corollary}
\begin{proof}
Such an algorithm, or such a code construction followed by one $V$-test, puts $\SAT$ in $\Pclass$. Proposition~\ref{prop:sat} would give $\Pclass=\NPclass$, contrary to Theorem~\ref{thm:pairing}.
\end{proof}

\begin{remark}[Role of the hypotheses]
The proof uses a binary address to locate an error, but its outer witness is an arbitrary structure element containing the complete trace. It does not prove a separation between deterministic and binary-nondeterministic polynomial time. It also uses unit-cost equality of arbitrary elements, notably of complete configuration codes. The theorem concerns this specified computation model, not a model in which the cost of comparing or guessing a tree is its expanded size. In the auxiliary expansion
$\mathcal N^{\mathrm{pair}}$ the trace witness is one natural number,
even when its usual binary representation is very long. The trace
inclusion of Proposition~\ref{prop:exp-mixed} applies to that expansion;
its unconditional transfer to the original $\{+,\cdot,0,1\}$-reduct
is not asserted. Corollary~\ref{cor:elementary-arithmetic} gives an
unconditional transfer in the elementary arithmetic signature.
The result is about unit-cost computation over the displayed structure
and does not settle the classical finite-alphabet
$\Pclass$ versus $\NPclass$ question.
\end{remark}

\section{Arithmetic, Malcev's algebra, and the absence of projections}
\label{sec:malcev}

\subsection{Simulating operations by their graphs}

We next remove $L,R$ from the hypotheses on the original language. We use decision-class simulations; no passage from a decidable graph to a deterministic witness-producing function is assumed.

\begin{lemma}[Simulation by graphs]
\label{lem:graphs}
Assume $\Pclass(\Sstr)=\NPclass(\Sstr)$. Let $\Sstr^+$ be obtained by adjoining finitely many total functions whose graphs belong to $\NPclass(\Sstr)$. Then, as decision classes on the same universe,
\[
 \NPclass(\Sstr^+)=\Pclass(\Sstr^+)=\Pclass(\Sstr).
\]
\end{lemma}
\begin{proof}
Under the assumed collapse, every added graph has a deterministic polynomial-time decider over $\Sstr$. Simulate a nondeterministic polynomial-time $\Sstr^+$-machine over $\Sstr$. For every call to an added function, guess its proposed output and check the graph. A failed check rejects that simulation branch. Totality ensures that a correct guess exists; single-valuedness ensures that every successful guess is the intended output. All later tests, including negative relation tests, are therefore performed on the correct data. Only polynomially many calls and checks occur. Thus
\[
 \NPclass(\Sstr^+)\subseteq\NPclass(\Sstr)
 =\Pclass(\Sstr)\subseteq\Pclass(\Sstr^+)
 \subseteq\NPclass(\Sstr^+).
\]
These inclusions prove the assertion.
\end{proof}

\begin{theorem}[Injective binary operations]
\label{thm:injective}
Let $\Sstr$ be a finite-signature structure with two distinct available elements and a jointly injective binary primitive $*:S^2\to S$. Then
\[
 \Pclass(\Sstr)\ne\NPclass(\Sstr).
\]
This holds under either constant convention and in every finite expansion retaining $*$.
\end{theorem}
\begin{proof}
Suppose that the classes collapse. The image relation
\[
 \operatorname{Node}(t)\quad\Longleftrightarrow\quad
 \exists u,v\;(t=u*v)
\]
belongs to $\NPclass(\Sstr)$ and hence to $\Pclass(\Sstr)$; its complement has a polynomial-time decider as well. Define total projections by taking the unique components on the image of $*$ and returning $t$ off that image. Their graphs are
\begin{align*}
 \Gamma_L(t,u)&\quad\Longleftrightarrow\quad
       (\exists v\;t=u*v)\lor(\neg\operatorname{Node}(t)\land u=t),\\
 \Gamma_R(t,v)&\quad\Longleftrightarrow\quad
       (\exists u\;t=u*v)\lor(\neg\operatorname{Node}(t)\land v=t).
\end{align*}
Both graphs are in $\NPclass(\Sstr)$: choose the relevant alternative, guess its component when needed, or use the deterministic complement-of-image test. Injectivity makes the selected components unique.

By Lemma~\ref{lem:graphs}, adjoining $L,R$ preserves the assumed collapse. Their definitions give $L(u*v)=u$ and $R(u*v)=v$, contradicting Theorem~\ref{thm:pairing}. The argument applies in every finite expansion. Lemma~\ref{lem:constants} then gives the parameter-allowing conclusion.
\end{proof}

\subsection{Arithmetic and truth predicates}
\label{subsec:arithmetic-truth}

\subsubsection{Standard arithmetic}

\begin{corollary}[Standard arithmetic and its finite expansions]
\label{cor:arithmetic}
\label{cor:cantor}
Every finite expansion $\Sstr$ of
$\mathcal N=(\mathbb N;+,\cdot,0,1)$ satisfies
$\Pclass(\Sstr)\ne\NPclass(\Sstr)$, under either convention for program
constants. In particular, for every $V\subseteq\mathbb N$,
\[
 \boxed{\Pclass(\mathcal N_V)\ne\NPclass(\mathcal N_V).}
\]
Consequently, no such $V$ admits a uniformly polynomial-time code
construction $C_\psi(\bar a)$ over $\mathcal N_V$ satisfying
\begin{equation}
 C_\psi(\bar a)\in V
 \quad\Longleftrightarrow\quad
 \mathcal N_V\models\psi(\bar a)
 \label{eq:arithmetic-fast-truth}
\end{equation}
for all existential formulas of the expanded language
$\{+,\cdot,0,1,V\}$ and all parameter tuples. The matrix may contain
negative occurrences of $V$, and the input to the encoder consists of
the explicit formula and the tuple.
\end{corollary}
\begin{proof}
Suppose that a finite expansion $\Sstr$ of $\mathcal N$ has a collapse.
Formula~\eqref{eq:cantor-graph} and the two projection formulas define
the graphs of $c,L_c,R_c$ by existential formulas in the original
arithmetic language. Their graphs therefore lie in $\NPclass(\Sstr)$.
By Lemma~\ref{lem:graphs}, adjoining these three total functions
preserves the assumed equality. The resulting structure satisfies
\eqref{eq:pairing}, contradicting Theorem~\ref{thm:pairing}.
This argument works for every finite expansion, so
Lemma~\ref{lem:constants} also excludes a collapse with arbitrary
finitely many program constants. Equivalently, the conclusion is an
instance of the first-order definable pairing obstruction stated in
Corollary~\ref{cor:definable-pairing} below.

Finally, a polynomial-time encoder satisfying
\eqref{eq:arithmetic-fast-truth}, followed by one $V$-test, would decide
$\SAT(\mathcal N_V)$ in polynomial time. Proposition~\ref{prop:sat}
would imply the excluded collapse. No arithmetical definability or
ordinary computability assumption on $V$ is used; the encoder itself
may use $V$ as a primitive test.
\end{proof}

There is also a direct term version of the arithmetic argument:
\[
 j(x,y)=(x+y)\cdot(x+y+1)+(y+y)=2c(x,y)
\]
is jointly injective and is computed by a fixed arithmetic term.
Adjoining $j$ changes neither polynomial-time class, since every call
is replaced by a bounded number of arithmetic instructions;
Theorem~\ref{thm:injective} therefore gives the same contradiction.
The pairing-only structures $(\mathbb N;c,0,1)$ and
$(\mathbb N;c,L_c,R_c,0,1)$ also satisfy the separation by the earlier
theorems, but the corollary above concerns the original arithmetic
signature and all its finite expansions.

For conventional effective syntax coding, the claim about existential
truth also has a direct diagonal form with no efficiency hypothesis.
This statement is separate from the more flexible coding condition
\eqref{eq:arithmetic-fast-truth}. We first distinguish a formula, its
numerical code, and a numeral denoting that code.

\begin{definition}[Syntax codes, numerals, and diagonal substitution]
\label{def:godel-substitution}
Let $\mathcal L$ be one of the arithmetic languages considered here,
including the unary relation symbol $V$, and fix an injective
G\"odel coding of its finite syntactic expressions by natural numbers.
For an expression $\alpha$, the notation $\ulcorner\alpha\urcorner$
denotes its numerical code. For $n\in\mathbb N$, fix the closed numeral
\[
 \underline0=0,\qquad
 \underline{n+1}=(\underline n+1).
\]
Thus $\underline n$ is a term, whereas $n$ is a number, and the value
of $\underline n$ in every arithmetic expansion considered here is
$n$. The notation $\ulcorner\alpha\urcorner$ is used in the
metalanguage; its code is named inside a formula by the corresponding
numeral.

Fix a distinguished variable $x$. A formula is admissible for the
following substitution operation if its free variables are contained
in $\{x\}$. If $q=\ulcorner\alpha\urcorner$ is the code of such a
formula, define
\[
 \operatorname{Sub}_{\mathcal L}(q,n)
       =\ulcorner\alpha[\underline n/x]\urcorner.
\]
Here $\alpha[\underline n/x]$ is obtained by replacing every free
occurrence of $x$ by the numeral $\underline n$; bound occurrences are
unchanged. The numeral is closed, so this substitution cannot cause
variable capture, and the resulting formula is a sentence.

If $q$ is not the code of an admissible formula, we define
$\operatorname{Sub}_{\mathcal L}(q,n)=0$ for every $n$. This default
output is a number, not a truth value; it only makes the numerical
function defined on all inputs. Define the diagonal-substitution
function by
\[
 d_{\mathcal L}(q)=\operatorname{Sub}_{\mathcal L}(q,q).
\]
In particular, if $q=\ulcorner\alpha(x)\urcorner$ is admissible, then
\[
 d_{\mathcal L}(q)=\ulcorner\alpha[\underline q/x]\urcorner.
\]
These operations manipulate finite syntax only. They never use the
interpretation of $V$.
\end{definition}

For an effective coding, expressions can be encoded and decoded
algorithmically, well-formedness and free-variable tests are decidable,
and numeral substitution is computable. The above convention then
makes $\operatorname{Sub}_{\mathcal L}$ and $d_{\mathcal L}$ total
computable numerical functions. We will use a stronger, elementary
version of these requirements in
Section~\ref{subsubsec:elementary-arithmetic}.

\begin{proposition}[No existential truth predicate for its own arithmetic language]
\label{prop:arithmetic-diagonal}
Fix an effective G\"odel coding of
$\mathcal L_{\mathrm{ar}}(V)=\{+,\cdot,0,1,V\}$ with the recognition
and substitution properties just specified. There is no
$V\subseteq\mathbb N$ such that
\begin{equation}
 \ulcorner\sigma\urcorner\in V
 \quad\Longleftrightarrow\quad
 \mathcal N_V\models\sigma
 \label{eq:arithmetic-self-truth}
\end{equation}
for every existential sentence $\sigma$ of this same expanded language,
allowing negative $V$-atoms in its quantifier-free matrix.
\end{proposition}
\begin{proof}
Write $\operatorname{Sub}_{\mathrm{ar}}$ and $d_{\mathrm{ar}}$ for the
functions of Definition~\ref{def:godel-substitution} in this language.
The admissible formulas here may have arbitrary quantifiers; only their
free variables are restricted. On input $q$, an ordinary algorithm
checks admissibility, constructs the numeral $\underline q$, and
substitutes it for the free occurrences of $x$, returning $0$ in the
inadmissible case. Thus $d_{\mathrm{ar}}$ is total computable, regardless
of the interpretation of $V$.

Its graph is a computably enumerable subset of $\mathbb N^2$. By the
Davis--Putnam--Robinson--Matiyasevich theorem \cite{Matiyasevich1971},
there is a quantifier-free formula $\delta(x,z,\bar u)$ in the
\emph{original} arithmetic language $\{+,\cdot,0,1\}$ such that, for
all $a,z\in\mathbb N$,
\begin{equation}
 \mathcal N\models\exists\bar u\,\delta(a,z,\bar u)
 \quad\Longleftrightarrow\quad z=d_{\mathrm{ar}}(a).
 \label{eq:diagonal-graph}
\end{equation}
Here $a,z$ on the left specify the values assigned to the two free
variables. For each $a$, witnesses $\bar u$ exist when
$z=d_{\mathrm{ar}}(a)$ and do not exist for any other $z$. The output
$z$ is unique; the auxiliary witnesses $\bar u$ need not be unique.

Choose the displayed variables distinct and form the formula
\[
 \beta(x)=\exists z\,\exists\bar u\,
       \bigl(\delta(x,z,\bar u)\land\neg V(z)\bigr).
\]
It has no free variable other than $x$ and is existential in
$\mathcal L_{\mathrm{ar}}(V)$. Since $\delta$ contains no $V$, its
interpretation is the same in $\mathcal N$ and $\mathcal N_V$.

\begin{samepage}
Equation~\eqref{eq:diagonal-graph} therefore gives, for every
$a\in\mathbb N$,
\[
\begin{aligned}
 \mathcal N_V\models\beta[\underline a/x]
 &\quad\Longleftrightarrow\quad
   \exists z\in\mathbb N\,
        \bigl(z=d_{\mathrm{ar}}(a)\ \text{and}\ z\notin V\bigr)\\
 &\quad\Longleftrightarrow\quad d_{\mathrm{ar}}(a)\notin V.
\end{aligned}
\]
\end{samepage}

Now take the numerical code of this already constructed formula and
then perform the indicated substitution:
\[
 e=\ulcorner\beta(x)\urcorner,\qquad
 \lambda=\beta[\underline e/x].
\]
The formula $\lambda$ is a closed existential sentence. By the
\emph{definition} of diagonal substitution,
\[
 d_{\mathrm{ar}}(e)
   =\operatorname{Sub}_{\mathrm{ar}}(e,e)
   =\ulcorner\beta[\underline e/x]\urcorner
   =\ulcorner\lambda\urcorner.
\]
Thus $e$ codes the one-variable formula, while $d_{\mathrm{ar}}(e)$
codes the sentence obtained by substituting the numeral for that
formula's code. Specializing the preceding semantic equivalence to
$a=e$ yields
\[
 \mathcal N_V\models\lambda
 \quad\Longleftrightarrow\quad
 \ulcorner\lambda\urcorner\notin V.
\]
If \eqref{eq:arithmetic-self-truth} held, applying it to this particular
existential sentence would give
\[
 \ulcorner\lambda\urcorner\in V
 \quad\Longleftrightarrow\quad
 \mathcal N_V\models\lambda
 \quad\Longleftrightarrow\quad
 \ulcorner\lambda\urcorner\notin V,
\]
a contradiction. There is no circular definition: the numerical
function $d_{\mathrm{ar}}$ and its arithmetic graph representation are
fixed first, then $\beta$, then its code $e$, and finally $\lambda$.
Diophantine representation is what keeps $\lambda$ in the purely
existential fragment.
\end{proof}

\begin{remark}[Comparison with Tarski and other truth predicates]
\label{rem:tarski}
Tarski's undefinability theorem excludes an arithmetic formula defining
the full truth set of the original arithmetic language
\cite{Tarski1956}; the distinction between object language and
metalanguage is also explained in \cite[Sections~7--10]{Tarski1944}.
It does not prevent the external set of true arithmetic sentences from
being named by a new predicate. Restricted truth predicates are
possible as well: the true existential sentences of the original
language form a computably enumerable set, hence an existentially
definable set by the same Diophantine representation theorem
\cite{Matiyasevich1971}. Such predicates describe truth for the original
language; they do not satisfy \eqref{eq:arithmetic-self-truth} for the
language containing the predicate itself. Negating an arithmetic
definition of existential truth need not preserve existential form,
whereas $\neg V(z)$ is quantifier-free in the expanded language.

Proposition~\ref{prop:arithmetic-diagonal} is a Tarski-style diagonal
obstruction for standard effective codes, with no complexity bound.
Corollary~\ref{cor:arithmetic} instead excludes every uniformly fast
code construction in the unit-cost model, even one using an arbitrary
primitive $V$ and therefore not necessarily computable in the ordinary
sense. The main pairing proof requires neither Diophantine
representation nor an arithmetization of syntactic self-reference.
Kripke's fixed-point truth predicates allow truth-value gaps
\cite{Kripke1975}; they do not provide the total classical predicate
required by the displayed biconditionals. The positive unary
construction in Section~\ref{sec:positive} uses a different structure
and coding scheme; this arithmetic diagonal argument does not by
itself apply to it.
\end{remark}

\subsubsection{Elementary arithmetic and term representations}
\label{subsubsec:elementary-arithmetic}

Let $\mathsf{Elem}$ denote the class of Kalm\'ar elementary functions:
the smallest class containing the constants, coordinate projections,
addition, and truncated subtraction, and closed under substitution,
bounded summation, and bounded product. Mazzanti's substitution-basis
results represent these functions by fixed arithmetic terms
\cite{Mazzanti2002}. We use the three-operation basis of Prunescu,
Sauras-Altuzarra and Shunia \cite[Theorem~1]{PrunescuSaurasShunia2026}:
\begin{equation}
 \mathsf{Elem}
   =\bigl\langle x+y,\ x\bmod y,\ 2^x\bigr\rangle_{\mathrm{sub}}.
 \label{eq:elementary-basis}
\end{equation}
Here the brackets mean closure under substitution, with constants and
coordinate projections available. Define
\[
 \mathcal N_{\mathrm{el}}
   =(\mathbb N;+,\operatorname{mod},\operatorname{exp}_2,0,1),
 \qquad \operatorname{exp}_2(x)=2^x,
\]
using the total convention $x\bmod0=x$. Every constant needed by a
fixed term can be built from $0,1$ and addition, so the basis theorem
does not require an infinite computational signature.

Cantor pairing $c$ and its two projections are elementary. For the
projections, the inverse description in
Section~\ref{subsec:running-examples} can be written as
\[
 s(n)=\sum_{k=1}^{n}\mathbf 1_{\{k(k+1)\leq2n\}},
 \qquad
 r(n)=n-\frac{s(n)(s(n)+1)}2,
\]
\[
 L_c(n)=s(n)-r(n),\qquad R_c(n)=r(n).
\]
The symbol $\mathbf 1$ denotes the characteristic function of the
indicated elementary relation. The bounded sum and the remaining
functions are elementary; all displayed subtractions are nonnegative.
Equation~\eqref{eq:elementary-basis} therefore supplies fixed terms
$t_c(x,y),t_L(z),t_R(z)$ in the language of $\mathcal N_{\mathrm{el}}$
representing $c,L_c,R_c$. In particular,
\[
 t_L(t_c(x,y))=x,\qquad t_R(t_c(x,y))=y.
\]
Each term has a fixed finite number of operation occurrences. Its
unit-cost evaluation is uniformly bounded, independently of the values
of its arguments. No uniform bound on the sizes of terms representing
\emph{all} elementary functions is needed.

This conclusion depends on the signature. Every unary term in
$\{+,\cdot,0,1\}$ is a polynomial with nonnegative integer coefficients
and is therefore nondecreasing. But
\[
 L_c(1)=1>L_c(2)=0,\qquad R_c(2)=1>R_c(3)=0.
\]
Thus the projections are not terms in that original language. Their
ordinary computability and first-order definability are compatible
with this fact.

\begin{corollary}[Elementary arithmetic and its finite expansions]
\label{cor:elementary-arithmetic}
For every finite expansion $\Sstr$ of $\mathcal N_{\mathrm{el}}$,
\begin{gather*}
 \Pclass(\Sstr)\subsetneq\EXPclass(\Sstr)
 \subseteq\Sigma_2^{S,B}(\Sstr)\subseteq\PHclass(\Sstr),\\
 \Pclass(\Sstr)\ne\NPclass(\Sstr).
\end{gather*}
In particular,
$\Pclass(\Sstr)\subsetneq\Sigma_2^{S,B}(\Sstr)$.
These statements hold for arbitrary interpretations of the additional
primitives, including an arbitrary unary predicate $V$. Either
convention for finitely many program constants gives the same classes
here.
\end{corollary}
\begin{proof}
Adjoin $c,L_c,R_c$ as primitive functions. Replace each call to an added
function by evaluation of its fixed defining term. This incurs only
constant overhead and preserves $\Pclass$, $\NPclass$, $\EXPclass$,
and all fixed levels of the full and mixed polynomial hierarchies.
For the hierarchies, the quantified domains and polynomial bounds are
unchanged, and the deterministic matrix computation is simulated by
the same replacement. These simulations require no collapse
hypothesis and remain valid with every additional primitive interpreted
as in $\Sstr$.

The expanded structure satisfies \eqref{eq:pairing}. Apply
Proposition~\ref{prop:exp-mixed} and Theorem~\ref{thm:pairing}, and
transfer their conclusions back by this simulation. The strict
inclusion $\Pclass\subsetneq\EXPclass$ is
Theorem~\ref{thm:hierarchy}. Finally, every extra program constant is a
natural number and can be initialized by a fixed numeral, so allowing
finitely many such constants does not change the classes.
\end{proof}

The same argument applies in any finite-signature structure with two
available distinct elements whenever pairing and its projections are
represented by fixed terms satisfying \eqref{eq:pairing}; uniformly
bounded-time implementations suffice as well. Multiplication itself
is a term function of $\mathcal N_{\mathrm{el}}$ by
\eqref{eq:elementary-basis}, so its separation also follows from
Corollary~\ref{cor:arithmetic} after adjoining multiplication. The
additional conclusion here is the unconditional simulation of the
projections and the resulting trace inclusion in the original
$\mathcal N_{\mathrm{el}}$ signature.

The elementary basis also sharpens the syntactic diagonal argument.
We use the distinction between codes and numerals from
Definition~\ref{def:godel-substitution}, now for the language
\[
 \mathcal L_{\mathrm{el}}(V)
   =\{+,\operatorname{mod},\operatorname{exp}_2,0,1,V\}.
\]
The required syntax operations will be \emph{elementary} numerical
functions, rather than merely computable ones.

Here is one concrete way to obtain such a coding. Choose a finite
alphabet for expressions, including digits for variable indices, and
represent its symbols by digits $a_i\in\{0,\ldots,B-1\}$ for a fixed
base $B\geq2$ large enough for the alphabet. Encode a string
$a_0\cdots a_{m-1}$ by
\[
 \operatorname{enc}(a_0\cdots a_{m-1})
   =B^m+\sum_{i=0}^{m-1}a_iB^{m-1-i}.
\]
The initial base-$B$ digit $1$ records the length, including any leading
zero digits of the string. Fully parenthesized expressions and explicit
variable indices give decidable syntax and free-variable tests.
Use the unary numerals of Definition~\ref{def:godel-substitution}.
If a formula has code $q$, its length is $O(\log(q+2))$, whereas the
numeral $\underline n$ has $O(n+1)$ symbols. Replacing the free
occurrences of $x$ therefore produces
$O((n+1)\log(q+2))$ symbols, with numerical code at most
\[
 B^{O((n+1)\log(q+2))}.
\]
Parsing, checking variable occurrences, and constructing the
substituted string can be carried out by ordinary algorithms whose
running time is bounded by a fixed polynomial in $q+n+2$. This is a
deliberately loose bound, since the input string has only
$O(\log(q+2))$ symbols. Their numerical outputs have elementary bounds
as well. Hence the syntax tests, numeral construction, and totalized
substitution functions in this coding are elementary. The following
proposition applies to any coding with these properties.

\begin{proposition}[No quantifier-free truth predicate in elementary arithmetic]
\label{prop:elementary-diagonal}
Put $\mathcal N_{\mathrm{el},V}=(\mathcal N_{\mathrm{el}},V)$ and fix a
G\"odel coding of $\mathcal L_{\mathrm{el}}(V)$ for which the recognition
of quantifier-free formulas whose free variables are contained in
$\{x\}$, and substitution of the numerals $\underline n$, are
elementary. Denote this restricted substitution function by
$\operatorname{Sub}_{\mathrm{qf}}(q,n)$, with the meaning given in
Definition~\ref{def:godel-substitution} on admissible formulas.

If $q$ does not code a quantifier-free formula with free variables
contained in $\{x\}$, define
$\operatorname{Sub}_{\mathrm{qf}}(q,n)=0$ for every $n$. This default
value only makes the numerical function defined on every input; it
does not assign a truth value to an invalid expression.

There is no $V\subseteq\mathbb N$ such that
\begin{equation}
 \ulcorner\sigma\urcorner\in V
 \quad\Longleftrightarrow\quad
 \mathcal N_{\mathrm{el},V}\models\sigma
 \label{eq:elementary-self-truth}
\end{equation}
for every quantifier-free sentence $\sigma$ of this same expanded
language, allowing negative $V$-atoms. No definability or ordinary
computability assumption is made on $V$.
\end{proposition}
\begin{proof}
First define the total numerical function
\[
 d_{\mathrm{qf}}(q)
   =\operatorname{Sub}_{\mathrm{qf}}(q,q).
\]
It is elementary by the hypothesis on substitution and closure under
composition. In particular, if $q=\ulcorner\alpha(x)\urcorner$ is
admissible, then
\[
 d_{\mathrm{qf}}(q)=\ulcorner\alpha[\underline q/x]\urcorner.
\]
This function depends only on the chosen coding and on finite syntax;
it never tests membership in $V$.

By \eqref{eq:elementary-basis}, choose a fixed term $t_d(x)$ using only
$+,\operatorname{mod},\operatorname{exp}_2,0,1$ such that
\begin{equation}
 t_d^{\mathcal N_{\mathrm{el}}}(a)=d_{\mathrm{qf}}(a)
 \qquad(a\in\mathbb N).
 \label{eq:diagonal-term}
\end{equation}
The left-hand side means the numerical value of this term when $x$ is
assigned $a$. The term is fixed before the next formula is formed and
contains no occurrence of $V$.

Form the quantifier-free formula
\[
 \beta(x)=\neg V(t_d(x)).
\]
When taking its code, write out the entire arithmetic term $t_d(x)$:
the notation $t_d$ is an abbreviation in the proof, not a new function
symbol in the coded language. Only after this expansion set
\[
 e=\ulcorner\beta(x)\urcorner,\qquad
 \lambda=\beta[\underline e/x]=\neg V(t_d(\underline e)).
\]
Because the numeral $\underline e$ is closed and $\beta$ has no other
free variables, $\lambda$ is a quantifier-free sentence of the same
language. Its code is
\[
 d_{\mathrm{qf}}(e)
   =\operatorname{Sub}_{\mathrm{qf}}(e,e)
   =\ulcorner\lambda\urcorner.
\]
The numeral $\underline e$ denotes the number $e$, so
\eqref{eq:diagonal-term} gives the following value for the closed term
occurring inside $\lambda$:
\[
 \bigl(t_d(\underline e)\bigr)^{\mathcal N_{\mathrm{el}}}
   =t_d^{\mathcal N_{\mathrm{el}}}(e)
   =d_{\mathrm{qf}}(e)
   =\ulcorner\lambda\urcorner.
\]
Thus the semantics of the single negative $V$-atom in $\lambda$ gives
\[
 \mathcal N_{\mathrm{el},V}\models\lambda
 \quad\Longleftrightarrow\quad
 \ulcorner\lambda\urcorner\notin V.
\]
Applying the proposed truth condition
\eqref{eq:elementary-self-truth} to $\lambda$ would instead require
\[
 \ulcorner\lambda\urcorner\in V
 \quad\Longleftrightarrow\quad
 \mathcal N_{\mathrm{el},V}\models\lambda
 \quad\Longleftrightarrow\quad
 \ulcorner\lambda\urcorner\notin V,
\]
a contradiction. The construction uses the substitution function first,
its representing term second, and only then the formula and its code.
The term representation replaces the existential graph witnesses used
in Proposition~\ref{prop:arithmetic-diagonal}, so the diagonal sentence
here remains quantifier-free.
\end{proof}

Proposition~\ref{prop:elementary-diagonal} concerns the specified
syntax coding. It does not exclude arbitrary $V$-dependent encoders of
quantifier-free truth: if fixed numbers $a\in V$ and $b\notin V$ are
available, such an encoder can evaluate a quantifier-free input formula
using $V$ and return $a$ or $b$ according to its truth value.
Corollary~\ref{cor:elementary-arithmetic}, together with
Proposition~\ref{prop:sat}, still excludes a uniformly polynomial-time
encoder of \emph{existential} truth in the full expanded language.
The new diagonal proposition strengthens the fragment in
Proposition~\ref{prop:arithmetic-diagonal} by changing the arithmetic
signature and imposing the elementary coding hypothesis; the
comparison with external and partial truth predicates in
Remark~\ref{rem:tarski} continues to apply.

\subsection{The negative answer to the final \emph{Les petits cailloux} question}

Recall the free Malcev algebra of Definition~\ref{def:malcev}, for an arbitrary nonempty generator set $G$. Its universe $T_G$ is defined recursively: every $g\in G$ is a tree, and for $u,v\in T_G$ the expression $u*v$ is a tree. Generators are distinct from all products, and two products are equal if and only if their respective left and right subtrees are equal. Every element is a finite term. This is the Malcev algebra considered here; it is not the anticommutative algebraic structure also called a Malcev algebra in other contexts.

Choose two distinct named elements $0,1\in T_G$ and write
\[
 \M=(T_G;*,0,1),\qquad
 \M_V=(T_G;*,V,0,1).
\]
When $G$ contains two designated generators they may be used as $0,1$. For a one-generator algebra, one may instead use that generator and its product with itself. The proof does not require the named elements to be leaves. Define total projections by
\[
 L(u*v)=u,\qquad R(u*v)=v,\qquad L(g)=R(g)=g\quad(g\in G),
\]
and put $\M_V^+=(T_G;*,L,R,V,0,1)$.

\begin{corollary}[Negative answer to Question~\ref{question:poizat}]
\label{cor:poizat}
For every unary predicate $V\subseteq T_G$,
\[
 \boxed{\Pclass(\M_V)\ne\NPclass(\M_V),\qquad
        \Pclass(\M_V^+)\ne\NPclass(\M_V^+).}
\]
Consequently there is no truth-predicate expansion of Malcev's algebra answering Poizat's free-tree question positively, in the unit-cost framework of \emph{Les petits cailloux}. The statement holds for every cardinality of $G$, in every finite further expansion, and under either constant convention.
\end{corollary}
\begin{proof}
The pairing function of $T_G$ is jointly injective, so Theorem~\ref{thm:injective} applies to $\M_V$. The expanded structure $\M_V^+$ satisfies \eqref{eq:pairing}, so Theorem~\ref{thm:pairing} applies to it directly. Both results permit arbitrary additional primitives in a fixed finite signature, including the chosen interpretation of $V$. The final assertions are part of those theorems.
\end{proof}

\begin{corollary}[The two projection conventions]
\label{cor:projection-equivalence}
For every $V\subseteq T_G$,
\[
 \bigl[\Pclass(\M_V)=\NPclass(\M_V)\bigr]
 \quad\Longleftrightarrow\quad
 \bigl[\Pclass(\M_V^+)=\NPclass(\M_V^+)\bigr].
\]
Both statements are false. The implication from a collapse in the pairing language to a collapse in the projection language is already supplied by a computation simulation.
\end{corollary}
\begin{proof}
Under a hypothetical collapse for $\M_V$, the complement of $\operatorname{Node}(t)$ is decidable in polynomial time. In the free algebra it is exactly the set of generators. Thus the two graphs in the proof of Theorem~\ref{thm:injective} are the graphs of the specified projections fixing generators, and Lemma~\ref{lem:graphs} gives a collapse for $\M_V^+$. Corollary~\ref{cor:poizat} rules out the latter and also rules out the former. This proves the displayed equivalence with both sides false.
\end{proof}

The equivalence is not an assertion that adding projections preserves the individual deterministic complexity classes without hypotheses. The simulation guesses projection outputs and uses the hypothetical collapse; it does not compute those outputs deterministically merely because their graphs are decidable.

\begin{corollary}[First-order definable pairing]
\label{cor:definable-pairing}
If a finite-signature structure with two distinct available elements satisfies $\Pclass=\NPclass$, there is no total jointly injective map $J:S^2\to S$ whose graph is first-order definable, even with finitely many parameters.
\end{corollary}
\begin{proof}
First name the finitely many definition parameters, preserving the collapse by Lemma~\ref{lem:constants}. By Proposition~\ref{prop:ph-collapse}, every fixed first-order definable relation is polynomial-time decidable under a collapse. Explicitly, induct on the fixed formula, using closure under Boolean combinations and the fact that existential projection of a polynomial-time relation is in $\NPclass=\Pclass$. Universal quantification follows by complementation. This is an induction for each fixed formula, not a uniform bound for arbitrary quantifier alternation.

The graph of $J$ is consequently in $\Pclass$ and hence in $\NPclass$. Lemma~\ref{lem:graphs} allows us to adjoin $J$ while retaining the collapse. Theorem~\ref{thm:injective} gives a contradiction. The parameter-allowing version follows by first applying Lemma~\ref{lem:constants} to a finite expansion with a fixed-signature collapse.
\end{proof}

% The positive constructions retain their full displayed signatures.

\section{Infinite structures with binary operations and \texorpdfstring{$\mathsf P=\mathsf{NP}$}{P = NP}}
\label{sec:positive}

We now recall the unary existence theorem and derive five structures with
binary operations.  In this section the letters $\mathcal A,\ldots,\mathcal F$
refer to these particular examples.  Each asserted class equality concerns
the entire displayed signature.

\begin{theorem}[Prunescu {\cite[Theorem~5.2]{PrunescuJSL}}]
\label{thm:unary-input}
There is a countably infinite structure
\[
  \mathcal A=(U;s_0,s_1,p,a,V)
\]
with three unary functions, one constant, and one unary predicate, such that
$\operatorname{SAT}(\mathcal A)$ is decidable in polynomial unit cost and
\[
  \mathsf P(\mathcal A)=\mathsf{NP}(\mathcal A).
\]
Its universe is a countably infinite disjoint union of rooted successor
trees: on each component $s_0$ and $s_1$ freely generate the finite words
from the root, $p(s_i(u))=u$, and $p$ fixes the root.  The constant $a$ names
one root.  In particular, $s_0,s_1$ are injective with disjoint images, $p$
is surjective, and all the values $s_0^j(a)$, $j\geq 0$, are distinct.
\end{theorem}

The theorem gives a quadratic bound in its explicit formula
encoding. Renaming the finitely many occurring variables to consecutive
indices converts the explicit finite-alphabet encodings used here to
that encoding with only polynomial overhead. Only the polynomial bound
and the displayed algebraic properties will be used here; the
construction of $V$ is not repeated.
The distinct ground-term values $s_0(a)$ and $s_1(a)$ can serve as the
binary values.  Naming finitely many fixed ground-term values does not
change any of the complexity classes.

\begin{lemma}[Replacing primitive instructions]
\label{lem:primitive-replacement}
Let two finite-signature structures have the same universe.  If every
primitive operation and relation of either structure is computable by a
fixed bounded-time program over the other, and the named constants can
also be initialized in bounded time, their deterministic and
nondeterministic polynomial-time classes coincide.
\end{lemma}
\begin{proof}
Replace each primitive instruction by its fixed simulation.  The running
time changes by at most a constant factor and a fixed initialization cost.
An arbitrary guessed element is the same witness in either structure, so
the same replacement also applies to nondeterministic computations.
\end{proof}

\subsection{A binary operation retaining the unary predicate}

Set
\[
  c_0=a,\qquad c_1=s_0(a),\qquad c_2=s_1(a),
\]
and define on $U$
\[
  x\star y=
  \begin{cases}
    s_0(y),&x=c_0,\\
    s_1(y),&x=c_1,\\
    p(y),&x=c_2,\\
    x,&x\notin\{c_0,c_1,c_2\}.
  \end{cases}
\]
Let $\mathcal B=(U;\star,V,a)$.

\begin{corollary}
\label{cor:example-B}
The countably infinite structure $\mathcal B$ has one binary operation,
one unary predicate, and one constant, and satisfies
$\mathsf P(\mathcal B)=\mathsf{NP}(\mathcal B)$.
Its binary operation depends essentially on both arguments, and the
subalgebra generated by $a$ is infinite.
\end{corollary}
\begin{proof}
The tags and the unary functions of $\mathcal A$ are recovered by
\[
 \begin{gathered}
 c_1=a\star a,\qquad c_2=(a\star a)\star a,\\
 s_0(y)=a\star y,\qquad s_1(y)=c_1\star y,\qquad p(y)=c_2\star y.
 \end{gathered}
\]
The predicate and the constant are unchanged.  Conversely, the displayed
case definition computes $\star$ over $\mathcal A$ using finitely many
equality tests and one old operation.  Thus Lemma~\ref{lem:primitive-replacement}
and Theorem~\ref{thm:unary-input} prove the class equality.

For a fixed first argument $a$, the map $y\mapsto a\star y=s_0(y)$ is
injective; hence the second argument is essential.  Two distinct first
arguments outside the finite tag set give distinct outputs for every
fixed second argument; hence the first argument is essential.  Iterating
the recovered $s_0$ on $a$ gives infinitely many elements of the generated
subalgebra.
\end{proof}

\subsection{One binary operation and one constant}

In addition to $c_0,c_1,c_2$, put $c_3=s_0^2(a)$.  These four elements are
distinct.  Define
\[
  x\diamond y=
  \begin{cases}
    s_0(y),&x=c_0,\\
    s_1(y),&x=c_1,\\
    p(y),&x=c_2,\\
    c_2,&x=c_3\text{ and }V(y),\\
    c_1,&x=c_3\text{ and }\neg V(y),\\
    x,&x\notin\{c_0,c_1,c_2,c_3\}.
  \end{cases}
\]
Let $\mathcal C=(U;\diamond,a)$; in particular, $V$ is not a primitive
predicate of $\mathcal C$.

\begin{corollary}
\label{cor:example-C}
There is a countably infinite algebra $\mathcal C$ with one binary
operation and one constant such that
$\mathsf P(\mathcal C)=\mathsf{NP}(\mathcal C)$.
Its operation depends essentially on both arguments, and its constant
generates an infinite subalgebra.
\end{corollary}
\begin{proof}
The unnamed tags are fixed ground-term values:
\[
 c_1=a\diamond a,\qquad
 c_2=(a\diamond a)\diamond a,\qquad
 c_3=a\diamond(a\diamond a).
\]
We recover all the primitives of $\mathcal A$ by
\[
 \begin{aligned}
 s_0(y)&=a\diamond y,&s_1(y)&=c_1\diamond y,&p(y)&=c_2\diamond y,\\
 V(y)&\quad\Longleftrightarrow\quad c_3\diamond y=c_2.
 \end{aligned}
\]
Each operation or test therefore requires a bounded number of
$\mathcal C$-instructions.  Conversely, the case definition of $\diamond$
uses only a bounded number of $\mathcal A$-instructions.  Apply
Lemma~\ref{lem:primitive-replacement} and
Theorem~\ref{thm:unary-input}.  The row with first argument $a$ is injective
in the second argument, whereas distinct first arguments outside the
four-element tag set give distinct outputs for any second argument.
Finally, the recovered iterates $s_0^j(a)$ belong to the subalgebra
generated by $a$ and are pairwise distinct.
\end{proof}

Neither $\star$ nor $\diamond$ is associative.  For example,
\[
 (a\diamond a)\diamond u=s_1(u)
 \ne s_0^2(u)=a\diamond(a\diamond u),
\]
using the disjoint images of $s_0$ and $s_1$; the same calculation applies
to $\star$.

\subsection{An associative rectangular-band operation}

On $U^2$ define
\[
 (u,v)\odot(w,z)=(u,z),\qquad J(u,v)=(v,u),
\]
and, for $f\in\{s_0,s_1,p\}$, define
\[
 \widehat f(u,v)=(f(u),f(v)).
\]
Put $W(u,v)\Longleftrightarrow V(u)$ and $a_{\mathcal D}=(a,a)$, and let
\[
 \mathcal D=
 (U^2;\odot,J,\widehat s_0,\widehat s_1,\widehat p,W,a_{\mathcal D}).
\]

\begin{corollary}
\label{cor:example-D}
The countably infinite structure $\mathcal D$ satisfies
$\mathsf P(\mathcal D)=\mathsf{NP}(\mathcal D)$.
Its binary operation is associative and idempotent, and satisfies
\[
 x\odot y\odot x=x,\qquad x\odot y\odot z=x\odot z.
\]
Every row and every column of its operation table has infinite image.
\end{corollary}
\begin{proof}
Let $L\in\mathsf{NP}(\mathcal D)$.  Replace every pair-valued data cell by
two $U$-valued cells.  The operation $\odot$ copies the first coordinate
of its first argument and the second coordinate of its second argument;
$J$ swaps two cells; each $\widehat f$ makes two calls to $f$; and $W$
queries $V$ on the first coordinate.  Pair equality is two equality
tests.  An arbitrary guessed pair is simulated by two arbitrary
$U$-element guesses.  This gives a polynomial-time nondeterministic
$\mathcal A$-machine for the flattened language
\[
 \widetilde L=
 \{(u_1,v_1,\ldots,u_n,v_n):
       ((u_1,v_1),\ldots,(u_n,v_n))\in L\}.
\]
Inputs of odd length are rejected in the flattened language.
Since $\mathsf P(\mathcal A)=\mathsf{NP}(\mathcal A)$, this language has a
deterministic polynomial-time $\mathcal A$-decider.

To run that decider inside $\mathcal D$, write $\Delta u=(u,u)$.
For every input $x=(u,v)$, the two diagonal values
\[
 x\odot J(x)=\Delta u,\qquad J(x)\odot x=\Delta v
\]
are computed internally in bounded time.  Thus the entire flattened
input can be prepared in diagonal form in linear time.  On diagonal
data, $\widehat f(\Delta u)=\Delta f(u)$,
$W(\Delta u)\Longleftrightarrow V(u)$, and equality is faithful.
The constant $a_{\mathcal D}=\Delta a$ supplies the original constant.
The deterministic $\mathcal A$-decider is therefore simulated in
polynomial time over $\mathcal D$.  Hence
$\mathsf{NP}(\mathcal D)\subseteq\mathsf P(\mathcal D)$, and the reverse
inclusion is automatic.

The algebraic assertions follow by direct coordinate calculation:
for $x=(u,v)$, $y=(w,z)$, and $t=(r,s)$, both bracketings of
$x\odot y\odot t$ equal $(u,s)$.  Also $x\odot x=x$ and
$x\odot y\odot x=x$.  For a fixed row $x=(u,v)$ its image is
$\{u\}\times U$, and for a fixed column $y=(w,z)$ its image is
$U\times\{z\}$.  Both are infinite.
\end{proof}

Although $\Delta u\odot\Delta v=(u,v)$, this operation does not encode
pairs of arbitrary elements of the full universe $U^2$ with recoverable
components.  Its recovery property is restricted to the diagonal subset.
The simulation uses the lifted unary operations in the displayed
signature.

\subsection{A commutative associative operation}

On $U$ define the flat meet with bottom $a$ by
\[
 x\wedge y=
 \begin{cases}
 x,&x=y,\\
 a,&x\ne y.
 \end{cases}
\]
Let $\mathcal E=(U;s_0,s_1,p,\wedge,a,V)$.

\begin{corollary}
\label{cor:example-E}
The countably infinite structure $\mathcal E$ satisfies
$\mathsf P(\mathcal E)=\mathsf{NP}(\mathcal E)$ and contains an
associative, commutative, idempotent binary operation depending
essentially on both arguments.
\end{corollary}
\begin{proof}
The new operation is computed over $\mathcal A$ by one equality test and
a copy or constant instruction, while every original primitive remains
available.  Thus Lemma~\ref{lem:primitive-replacement} and
Theorem~\ref{thm:unary-input} give the class equality.

Commutativity and idempotence follow immediately from the definition.
For any three arguments, both bracketings return their common value if
all three arguments are equal, and return $a$ otherwise.  This proves
associativity.  If $b\ne a$, then $b\wedge b=b$ whereas $a\wedge b=a$
and $b\wedge a=a$, which proves essential dependence on each argument.
\end{proof}

This corollary concerns the displayed expansion, which retains the unary
functions and the predicate of $\mathcal A$.

\subsection{One binary operation with all rows and columns infinite}

We finally compress the full signature of $\mathcal D$.  Put
\[
 F_1=\widehat s_0,\qquad F_2=\widehat s_1,\qquad
 F_3=\widehat p,\qquad F_4=J.
\]
For $j\geq0$, write $r_j=s_0^j(a)$, and choose the six distinct diagonal
tags and the auxiliary value
\[
 \begin{gathered}
 t_i=\Delta r_i\quad(1\leq i\leq4),\qquad
 t_E=\Delta r_5,\qquad t_V=\Delta r_6,\\
 k=r_7,\qquad \kappa=\Delta k.
 \end{gathered}
\]
These are fixed ground-term values in $\mathcal D$, and $k$ is different
from the first coordinate of every tag.  Define
\[
 E(u,v)=(k,v)=\kappa\odot(u,v),\qquad
 G(x)=
 \begin{cases}
 F_2(x),&W(x),\\
 F_1(x),&\neg W(x),
 \end{cases}
\]
and set
\[
 x\circ y=
 \begin{cases}
 F_i(x),&y=t_i\quad(1\leq i\leq4),\\
 E(x),&y=t_E,\\
 G(x),&y=t_V,\\
 x\odot y,&y\notin\{t_1,t_2,t_3,t_4,t_E,t_V\}.
 \end{cases}
\]
Let
\[
 \mathcal F=
 (U^2;\circ,a_{\mathcal D},t_1,t_2,t_3,t_4,t_E,t_V,\kappa).
\]
The final signature contains one binary operation and eight constants,
with no unary function or primitive predicate.

\begin{corollary}
\label{cor:example-F}
The countably infinite algebra $\mathcal F$ satisfies
$\mathsf P(\mathcal F)=\mathsf{NP}(\mathcal F)$.
Every row and every column of its binary operation has infinite image,
and its named constants generate an infinite subalgebra.
\end{corollary}
\begin{proof}
The case definition computes $\circ$ in bounded time over $\mathcal D$,
and all named constants of $\mathcal F$ are fixed ground-term values
there.  Conversely, the unary operations are recovered by
\[
 F_i(x)=x\circ t_i,\qquad E(x)=x\circ t_E.
\]
Because $E(y)$ has first coordinate $k$, it cannot equal any of the six
tags.  Therefore the original product is recovered by
\[
 x\odot y=x\circ(y\circ t_E).
\]
Indeed, the right-hand side is $x\odot E(y)=(x_1,y_2)$.
The predicate is recovered by
\[
 W(x)\quad\Longleftrightarrow\quad x\circ t_V=x\circ t_2.
\]
For the false case, $F_1(x)\ne F_2(x)$ follows from the disjoint images
of $s_0$ and $s_1$.  Thus every primitive of $\mathcal D$ is simulated
in bounded time.  Lemma~\ref{lem:primitive-replacement} and
Corollary~\ref{cor:example-D} imply the asserted class equality.

For fixed $x$, all arguments $y=(k,d)$, $d\in U$, are outside the tag
set and give the distinct outputs $(x_1,d)$; hence each row has infinite
image.  An ordinary column $y$ has image $U\times\{y_2\}$.
Each special column also has infinite image: $F_1,F_2$ are injective,
$F_3$ is surjective, $F_4$ is bijective, and $E$ has image
$\{k\}\times U$.  Finally, $G$ is injective, since each of its branches
is injective and the two branch images are disjoint.  This deals with
all columns.  Iterating the recovered $F_1$ on $a_{\mathcal D}$ gives
infinitely many pairwise distinct ground-term values.
\end{proof}

The operation $\circ$ is not associative.  Choose $u\ne v$, let
$x=(u,u)$ and $y=(u,v)$, and use the interchange tag $t_4$.  Both $y$
and $J(y)$ are off the diagonal and therefore are not tags.  Consequently,
\[
 (x\circ y)\circ t_4=(v,u),\qquad
 x\circ(y\circ t_4)=(u,u).
\]
These examples show that essential dependence on both arguments, infinite
rows and columns, and infinite generation by the named constants are
compatible with $\mathsf P=\mathsf{NP}$.  The obstruction proved earlier
requires a way of recovering both arbitrary input elements from one
encoded element. By Proposition~\ref{prop:ph-collapse}, all five positive
examples also satisfy $\PHclass=\Pclass$, and every fixed level of their
mixed hierarchies equals $\Pclass$.

\section{Concluding observations}
\label{sec:conclusion}

The negative and positive results identify a precise distinction. In Malcev's algebra the binary pairing function retains both arbitrary arguments, and the projections make them available at unit cost. Theorem~\ref{thm:injective} shows that omitting projections does not rescue a possible collapse. By contrast, the operations of the positive examples either dispatch among unary instructions or discard information from their inputs. Essential dependence on both arguments, associativity in an expanded signature, and infinite row and column images do not by themselves supply a recoverable pairing of the full universe.

The five constructions all inherit their collapse from the unary structure. No independent forcing construction is needed. Their displayed signatures are part of their statements: the rectangular-band and flat-meet examples retain specified unary operations, while the two compressed algebras $\mathcal C$ and $\mathcal F$ retain no primitive predicates and have nonassociative binary operations.

Theorem~\ref{thm:pairing} applies to every final interpretation of the additional primitives. It therefore rules out a truth-predicate construction on the free binary-tree algebra regardless of whether that predicate is defined by recursion, a genericity argument, or another method. The contradiction uses finite traces only and has no dependence on
countability or well-foundedness beyond what is needed to describe the
particular free-tree example. Standard arithmetic
$(\mathbb N;+,\cdot,0,1)$ illustrates the obstruction through its
first-order definable Cantor pairing and projections. The trace
argument applies to the auxiliary definitional expansion, and the
graph simulation transfers the contradiction to the original signature
and every finite expansion, including an arbitrary predicate $V$.
For structures with primitive pairing and projections, or fixed terms
representing them, polynomial-time access to a packed trace places
deterministic exponential time in
$\Sigma_2^{S,B}$, whereas a hypothetical $\Pclass=\NPclass$ would
collapse that mixed class, and the entire polynomial hierarchy, to
$\Pclass$. The comparison with Tarski separates this computational
obstruction from both external truth for the base language and the
direct diagonal obstruction to existential truth in the expanded
arithmetic language. In the elementary arithmetic signature
$\{+,\operatorname{mod},\operatorname{exp}_2,0,1\}$, the substitution-basis theorem
provides fixed terms for pairing and its projections, so the trace
inclusion holds without a collapse hypothesis. With elementary syntax
coding, the same term representation of diagonal substitution yields
the quantifier-free truth obstruction of
Proposition~\ref{prop:elementary-diagonal}.

\section*{Acknowledgments}

This paper was produced in an AI-assisted workflow. AI assistants helped explore candidate constructions, draft proofs, and write software for numerical experiments.

\end{document}